\documentclass{amsart}
\usepackage{a4wide}

\usepackage{amsmath,amsthm,amssymb}     
\usepackage{graphicx}     
\usepackage{hyperref} 
\usepackage{url}
\usepackage{amsfonts} 
\usepackage{bm}   

\def\<{\langle}
\def\>{\rangle}

\usepackage{multicol}

\usepackage{tikz}
\usetikzlibrary{arrows}
\usetikzlibrary{decorations.markings}
\usetikzlibrary{decorations.pathreplacing}
\usetikzlibrary{patterns}
\usetikzlibrary{shapes.geometric}
\usepackage{tikz-3dplot}
\usepackage{tkz-graph}
\usepackage{tikz-cd}

\DeclareMathOperator\PSL{PSL}

\DeclareMathOperator\Dic{Dic}
\DeclareMathOperator{\Alt}{A}
\DeclareMathOperator{\Cyc}{C}
\DeclareMathOperator{\Dih}{D}
\DeclareMathOperator{\Sym}{S}
\DeclareMathOperator{\Quat}{Q}
\DeclareMathOperator{\SA}{SA}
\DeclareMathOperator{\SD}{SD}
\DeclareMathOperator{\OD}{OD}

\DeclareMathOperator{\DQ}{DQ}
\DeclareMathOperator{\QD}{QD}

\DeclareMathOperator\Aut{Aut}

\DeclareMathOperator\cl{cl}
\DeclareMathOperator\Gal{Gal}

\def\C{\mathbb{C}}

\def\Z{\mathbb{Z}}
\def\0{\mathbf{0}}
\def\1{\mathbf{1}}

\usepackage{xcolor}
\usepackage{color}
\definecolor{darkblue}{rgb}{0, 0, .6}
\definecolor{grey}{rgb}{.7, .7, .7}
\colorlet{lightred}{red!30!white}
\colorlet{lightblue}{blue!30!white}
\colorlet{lightyellow}{yellow!30!white}
\colorlet{keylime}{green!10!white}
\colorlet{darkgreen}{green!50!black}
\colorlet{midgreen}{green!70!black}
\colorlet{lightgreen}{darkgreen!30!white}
\colorlet{lightpurple}{violet!60!white}
\colorlet{darkgrey}{black!70}
\definecolor{Red}{rgb}{.9,0,0}
\definecolor{Blue}{rgb}{0,0,.9}
\definecolor{Green}{HTML}{7EC636}
\definecolor{Purple}{HTML}{D287FF} 

\colorlet{lred}{red!50!white}
\definecolor{lblue}{rgb}{0.60, 0.75, 1}
\colorlet{lgreen}{green!40!white}
\colorlet{lyellow}{yellow!30!white}
\definecolor{lpurple}{rgb}{0.7 0.45 0.9}
\colorlet{lorange}{orange!50!white}
\definecolor{lightgrey}{rgb}{.85, .85, .85}
\definecolor{amethyst}{rgb}{0.54 0.17 0.89}
\definecolor{amber}{rgb}{1.0 0.49 0.0}

\colorlet{lightestred}{red!15!white}
\colorlet{lightestblue}{blue!15!white}
\colorlet{lightestgreen}{darkgreen!15!white}
\colorlet{lightestpurple}{amethyst!20!white}
\definecolor{faded}{rgb}{.75, .75, .75}
\definecolor{midgrey}{rgb}{.35, .35, .35}

\definecolor{lpurple}{rgb}{0.7 0.45 0.9}
\colorlet{tred}{red!70!white}
\colorlet{tblue}{blue!60!white}     
\colorlet{tgreen}{green!75!black}   
\colorlet{tyellow}{yellow!60!white}
\colorlet{tlime}{green!80!white}
\colorlet{tpink}{pink}
\colorlet{tpurple}{lpurple}
\colorlet{torange}{orange}
\colorlet{darkgreen}{green!50!black}
\definecolor{lightgrey}{rgb}{.85, .85, .85}
\definecolor{Red}{rgb}{.9,0,0}
\definecolor{Blue}{rgb}{0,0,.9}
\definecolor{Purple}{rgb}{0.54 0.17 0.89}
\definecolor{faded}{rgb}{.75, .75, .75}
\colorlet{rfaded}{red!30!white}
\colorlet{bfaded}{blue!30!white}
\definecolor{purple}{rgb}{0.54 0.17 0.89}

\definecolor{lpurple}{rgb}{0.7 0.45 0.9}
\colorlet{lorange}{orange!50!white}

\colorlet{lightestblue}{blue!15!white}
\colorlet{lightestred}{red!15!white}

\definecolor{conj-purple}{HTML}{E4B7FF}
\definecolor{conj-orange}{HTML}{FFD4B7}
\definecolor{Green}{HTML}{7EC636}
\colorlet{Purple}{amethyst}

\newcommand{\Alert}[1]{\textcolor{Red}{#1}}
\newcommand{\Balert}[1]{\textcolor{Blue}{#1}}
\newcommand{\Galert}[1]{\textcolor{darkgreen}{#1}}
\newcommand{\Palert}[1]{\textcolor{purple}{#1}}

\tikzstyle{v} = [circle, draw, fill=lightgrey,inner sep=0pt, 
  minimum size=6mm]
\tikzstyle{r} = [draw, very thick, Red, -stealth]        
\tikzstyle{b} = [draw, very thick, Blue, -stealth]       
\tikzstyle{g} = [draw, very thick, darkgreen, -stealth]  
\tikzstyle{rr} = [draw, very thick, Red]                 
\tikzstyle{bb} = [draw, very thick, Blue]                
\tikzstyle{gg} = [draw, very thick, darkgreen]           
\colorlet{f}{faded}

\theoremstyle{definition}

\begin{document}

\title{Dihedralizing the quaternions}

\subjclass[2020]{Primary 20D99; Secondary 05C25}
\keywords{Cayley graph, dicyclic, diquaternion, quaternion, semidihedral, semiabelian, subgroup lattice}

\author{Matthew Macauley}
\address{School of Mathematical and Statistical Sciences, Clemson University, Clemson, SC 29634} 
\email{macaule@clemson.edu}

\maketitle

\begin{abstract}
In this paper, we take the classic dihedral and quaternion groups and explore questions like  ``\emph{what if we replace $i=e^{2\pi i/4}$ in $\Quat_8$ with a larger root of unity?}'' and ``\emph{what if we add a reflection to $\Quat_8$?}''  The delightful answers reveal lesser-known families like the dicyclic, diquaternion, semidihedral, and semiabelian groups, which come to life with visuals such as Cayley graphs, cycle graphs, and subgroup lattices.
\end{abstract}

\section{A mysterious group of order 16}

\begin{figure}[!ht]
  \begin{tikzpicture}[node distance=1cm,shorten >= -2pt, shorten <= -2pt,scale=.95]
    \begin{scope}[shift={(0,0)}]
      \tikzstyle{every node}=[font=\small]
      \node(G) at (3.5,6) {\!\!\!\!\!\!\!\!\!\!\!\Palert{$G=\<r,s\>$}};
      \node(rs) at (3.5,4.5) {$\<rs\>$};
      \node(r2-s) at (1.75,4.5) {\Palert{$\<r^2,s\>$}};
      \node(r) at (5.25,4.5) {$\<r\>$};
      \node(r4-s) at (0,3) {\Palert{$\<r^4,s\>$}};
      \node(r2s) at (1.75,3) {\Palert{$\<r^2s\>$}};
      \node(r2) at (3.5,3) {\Palert{$\<r^2\>$}};
      \node(s) at (-1.75,1.5) {$\<s\>$};
      \node(r4s) at (0,1.5) {$\<r^4s\>$};
      \node(r4) at (1.75,1.5) {\Palert{$\<r^4\>$}};
      \node(1) at (0,0) {\Palert{$\<1\>$}};
      \draw[f] (1) to (s);
      \draw[f] (1) to (r4s);
      \draw[f] (1) to (r4);
      \draw[f] (s) to (r4-s);
      \draw[f] (r4s) to (r4-s);
      \draw[f] (r4) to (r4-s);
      \draw[f] (r4) to (r2s);
      \draw[f] (r4) to (r2);
      \draw[f] (r4-s) to (r2-s);
      \draw[f] (r2s) to (r2-s);
      \draw[f] (r2) to (r2-s);
      \draw[f] (r2) to (r);
      \draw[f] (r2) to (rs);
      \draw[f] (r2-s) to (G);
      \draw[f] (r) to (G);
      \draw[f] (rs) to (G);
    \end{scope}
    \begin{scope}[shift={(7.5,0)}]
      \tikzstyle{every node}=[font=\footnotesize]
      \node(G) at (3.5,6) {\hspace{2mm}\Palert{$xGx^{-1}$}};
      \node(rs) at (3.5,4.5) {$x\<rs\>x^{-1}$};
      \node(r2-s) at (1.75,4.5) {\Palert{$x\<r^2,s\>x^{-1}$}};
      \node(r) at (5.25,4.5) {$x\<r\>x^{-1}$};
      \node(r4-s) at (0,3) {\!\!\!\!\Palert{$x\<r^4,s\>x^{-1}$}};
      \node(r2s) at (1.75,3) {\hspace{3mm}\Palert{$x\<r^2\!s\>x^{-1}$}};
      \node(r2) at (3.5,3) {\hspace{5mm}\Palert{$x\<r^2\>x^{-1}$}};
      \node(s) at (-1.75,1.5) {$x\<s\>x^{-1}$};
      \node(r4s) at (0,1.5) {$x\<r^4s\>x^{-1}$};
      \node(r4) at (1.75,1.5) {\hspace{3mm}\Palert{$x\<r^4\>x^{-1}$}};
      \node(1) at (0,0) {\Palert{$x\<1\>x^{-1}$}};
      \draw[f] (1) to (s);
      \draw[f] (1) to (r4s);
      \draw[f] (1) to (r4);
      \draw[f] (s) to (r4-s);
      \draw[f] (r4s) to (r4-s);
      \draw[f] (r4) to (r4-s);
      \draw[f] (r4) to (r2s);
      \draw[f] (r4) to (r2);
      \draw[f] (r4-s) to (r2-s);
      \draw[f] (r2s) to (r2-s);
      \draw[f] (r2) to (r2-s);
      \draw[f] (r2) to (r);
      \draw[f] (r2) to (rs);
      \draw[f] (r2-s) to (G);
      \draw[f] (r) to (G);
      \draw[f] (rs) to (G);
    \end{scope}
  \end{tikzpicture}
  \caption{Conjugating the subgroup lattice of a mysterious group $G=\<r,s\>$ of order $16$. All subgroups but $\<s\>$ and $\<r^4s\>$ are either unicorns or have index $2$, and thus are normal by inspection. }\label{fig:mystery-group}
  \end{figure}
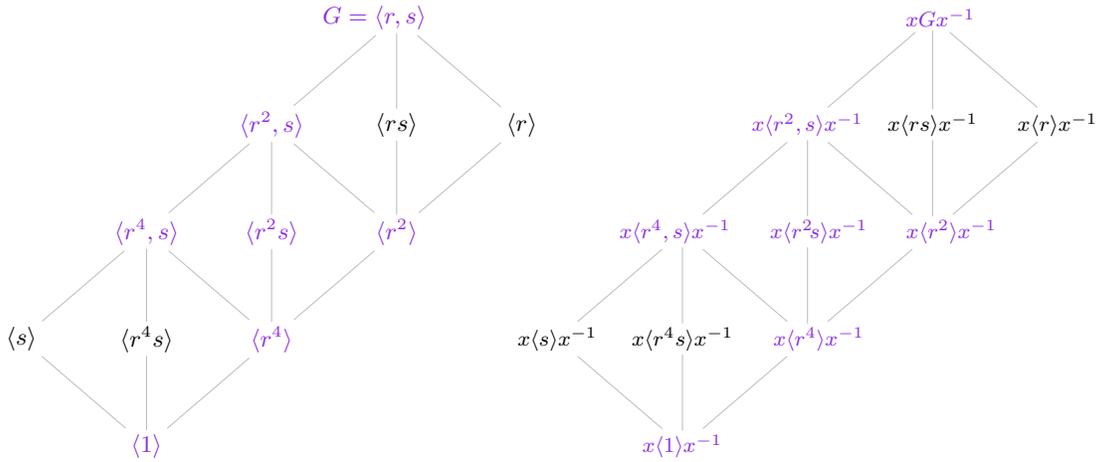

  We begin this article with a fun group theory puzzle. Figure~\ref{fig:mystery-group} (left) shows the subgroup lattice of a ``mystery group'' $G=\<r,s\>$ of order $16$, and we will see what we can deduce about it just by inspection. First, all three index-$2$ subgroups must be normal. Next, if we conjugate $G$ by some element $x\in G$, and compare the lattices in Figure~\ref{fig:mystery-group}, we see that any ``lattice automorphism invariant subgroup'' must be normal. For convenience, we will call such subgroups ``\emph{unicorns}.''

We have deduced that at least $9$ of the $11$ subgroups of $G$ must be normal, and the jury is still out on $\<s\>$ and $\<r^4s\>$. If the order-$2$ subgroup $\<s\>$ is normal, then it must be central, and $rs=sr$ forces $G$ to be abelian. Otherwise, it lies in a non-singleton conjugacy class, and the only other possible conjugate is $\<r^4s\>$. In this case, $\cl_G(\<s\>)=\big\{\<s\>,\<r^4s\>\big\}$, and $G$ is nonabelian. At this point, both of these remain in the realm of possibilities, and we challenge the reader, as a fun exercise, to find a group of order $16$ that has this subgroup lattice.

In the remainder of this article, we will explore new families of groups created from two classic examples: the dihedral group $\Dih_n$ and the quaternion group $\Quat_8$. None of these new families are particularly well known, primarily because their presentations are uninspiring. A \emph{Cayley graph} is a visualization of a presentation that makes multiplying elements very easy. With the right generating sets, this graph can bring these groups to life in new ways, highlight important structural features, and reveal surprising patterns. All groups that we will encounter in this article have very satisfying representations by $2\times 2$ matrices, which we will compare and contrast. By the end, we will discover the surprising answer to our mystery subgroup lattice in Figure~\ref{fig:mystery-group}.

\section{Quaternion and dihedral groups}

Let's begin by reviewing two classic examples of groups that appear in every algebra class. The \emph{quaternion group} $\Quat_8$ is defined as
\begin{equation}\label{eq:Q8-representation}
    \Quat_8=\big\<i,j\mid i^4=j^4=1,\,i^2=j^2=-1\big\>=\big\{\pm\!1,\pm i,\pm j,\pm k\big\}.
\end{equation}
Multiplication is anti-commutative, and $ij=k$. For the generalizations that we are about to introduce, it is helpful think of this group as being generated by: 
  \begin{itemize}
  \item a \Alert{$4^{\rm th}$ root of unity}, $i=\zeta_4=e^{2\pi i/4}$ (a $2\pi/4$-rotation in $\C$),
  \item the ``\Balert{imaginary number}'' $j$ from $\Quat_8$.
  \end{itemize}
  The standard $2\times 2$ complex representation of the quaternions is
  \begin{equation}\label{Q8-representation}
  \Quat_8\!=\!\big\<\Alert{i},\Balert{j},k\big\>\cong\left\<
  \underbrace{\Alert{\begin{bmatrix} i & 0 \\ 0 & -i\end{bmatrix}}}_{R=R_4},
  \underbrace{\Balert{\begin{bmatrix} 0 & -1 \\ 1 & 0\end{bmatrix}}}_{S},
  \underbrace{\begin{bmatrix} 0 & -i \\ -i & 0\end{bmatrix}}_{T=RS}\right\>.
  \end{equation}
  The \emph{dihedral group} $\Dih_n$ consists of the $2n$ symmetries of a regular $n$-gon. A standard generating set is 
  \[
  \Dih_n=\<r,f\>=\big\{\underbrace{1,r,r^2,\dots,r^{n-1}}_{\text{rotations}},\underbrace{f,rf,r^2f,\dots,r^{n-1}f}_{\text{reflections}}\big\},
  \]
  where $r$ is a $2\pi/n$ (counterclockwise) rotation, and $f$ is a reflection. We can think of $\Dih_n$ as being generated by:
  \begin{itemize}
  \item an \Alert{$n^{\rm th}$ root of unity}, $r=\zeta_n=e^{2\pi i/n}$ (a $2\pi/n$-rotation in $\C$),
  \item a \Balert{reflection} $f$.
  \end{itemize}
  This too has a canonical $2\times 2$ representation over $\C$: 
  \begin{equation}\label{eq:Dn-representation}
  \Dih_n\!=\!\<\Alert{r},\Balert{f}\>=\big\<r,f\mid r^n=f^2=1,\;rfr=f\big\>\cong\left\<\underbrace{\Alert{\begin{bmatrix}\zeta_n & 0
  \\ 0 & \bar{\zeta}_n\end{bmatrix}}}_{R_n},
  \Balert{\underbrace{\begin{bmatrix}0 & 1 \\ 1 & 0\end{bmatrix}}_{\Balert{F}}}\right\>\!.
  \end{equation}
  We call $F$ the ``complex reflection matrix'' because it maps $(z_1,z_2)\mapsto(z_2,z_1)$ in $\C^2$. 
  
  Throughout this article, we will exploit three effective but underutilized visual tools: Cayley graphs, cycle graphs, and subgroup lattices. All of these highlight features of groups, but in different ways. In a \emph{Cayley graph}, nodes represent group elements and edges represent generators. From each node are outgoing edges, one for each generator. Cayley graphs for $\Quat_8=\<i,j\>$ and $\Dih_6=\<r,f\>$ are shown in Figure~\ref{fig:Q8-D6-Cayley}, with $r$ denoting $i$ in $\Quat_8$. Throughout, we will write elements from \emph{left-to-right}, so arrows in a Cayley graph describe right-multiplication. Bidirected arrows will be drawn as undirected. When our groups contain elements that correspond to roots of unity, we will position these along the outer ring of the Cayley graph, and at the angle at which they are located on the unit circle. Though this online version of this article is in full color, this is done for added effects, and is largely unnecessary. In the few cases where grayscale figures would be potentially ambiguous, edges are dashed and/or dotted for extra emphasis in the (grayscale) printed journal version.

  \tikzstyle{R6-out} = [draw, very thick, Red,-stealth,bend right=18]
  \tikzstyle{R6-in} = [draw, very thick, Red,-stealth,bend left=12]
  \tikzstyle{R-out} = [draw, very thick, Red,-stealth,bend right=15]
  \tikzstyle{R-in} = [draw, very thick, Red,-stealth,bend left=12]
  \tikzstyle{B} = [draw, very thick, Blue,-stealth,bend right=25]
  \begin{figure}[!ht]
  \tikzstyle{v-small} = [circle, draw, fill=lightgrey,inner sep=0pt,
    minimum size=5.5mm]
  \begin{tikzpicture}[scale=1.2]
    \begin{scope}[shift={(0,0)},scale=1]
      \tikzstyle{every node}=[font=\footnotesize]
      \node (1) at (0:2) [v-small] {$1$};
      \node (r) at (90:2) [v-small] {$r$};
      \node (r2) at (180:2) [v-small] {$r^2$};
      \node (r3) at (270:2) [v-small] {$r^3$};
      \node (s) at (0:1) [v-small] {$j$};
      \node (sr) at (270:1) [v-small] {$r^3j$};
      \node (sr2) at (180:1) [v-small] {$r^2j$};
      \node (sr3) at (90:1) [v-small] {$rj$};
      \draw [b] (1) to (s); \draw [B] (s) to (r2);
      \draw [b] (r2) to (sr2); \draw [B] (sr2) to (1);
      \draw [b] (r) to (sr3); \draw [B] (sr3) to (r3);
      \draw [b] (r3) to (sr); \draw [B] (sr) to (r);
      \draw [r] (1) to [bend right] (r);
      \draw [r] (r) to [bend right] (r2);
      \draw [r] (r2) to [bend right] (r3);
      \draw [r] (r3) to [bend right] (1);
      \draw [r] (s) to [bend left=28] (sr);
      \draw [r] (sr) to [bend left=28] (sr2);
      \draw [r] (sr2) to [bend left=28] (sr3);
      \draw [r] (sr3) to [bend left=28] (s);
      \node at (0,0) {\normalsize $\Quat_8$};
    \end{scope}
    \begin{scope}[shift={(6,0)},scale=1]
      \tikzstyle{every node}=[font=\scriptsize]
      \node (1) at (0:2) [v-small] {$1$};
      \node (r) at (60:2) [v-small] {$r$};
      \node (r2) at (120:2) [v-small] {$r^2$};
      \node (r3) at (180:2) [v-small] {$r^3$};
      \node (r4) at (240:2) [v-small] {$r^4$};
      \node (r5) at (300:2) [v-small] {$r^5$};
      \node (s) at (0:1) [v-small] {$f$};
      \node (sr5) at (60:1) [v-small] {$rf$};
      \node (sr4) at (120:1) [v-small] {$r^2\!f$};
      \node (sr3) at (180:1) [v-small] {$r^3\!f$};
      \node (sr2) at (240:1) [v-small] {$r^4\!f$};
      \node (sr) at (300:1) [v-small] {$r^5\!f$};
      \path[bb] (1) to (s);
      \path[bb] (r) to (sr5);
      \path[bb] (r2) to (sr4);
      \path[bb] (r3) to (sr3);
      \path[bb] (r4) to (sr2);
      \path[bb] (r5) to (sr);
      \path[R-out] (1) to (r);
      \path[R6-out] (r) to (r2);
      \path[R6-out] (r2) to (r3);
      \path[R6-out] (r3) to (r4);
      \path[R6-out] (r4) to (r5);
      \path[R6-out] (r5) to (1);
      \path[R-in] (sr5) to (s);
      \path[R-in] (sr4) to (sr5);
      \path[R-in] (sr3) to (sr4);
      \path[R-in] (sr2) to (sr3);
      \path[R-in] (sr) to (sr2);
      \path[R-in] (s) to (sr);
      \node at (0,0) {\normalsize $\Dih_6$};
    \end{scope}
    \end{tikzpicture}
    \caption{Cayley graphs for the quaternion group $\Quat_8=\<i,j\>=\<r,j\>$ and the dihedral group $\Dih_6=\<r,f\>$. The blue arrows are undirected in $\Dih_6$ because $f$ has order $2$.}\label{fig:Q8-D6-Cayley}
  \end{figure}
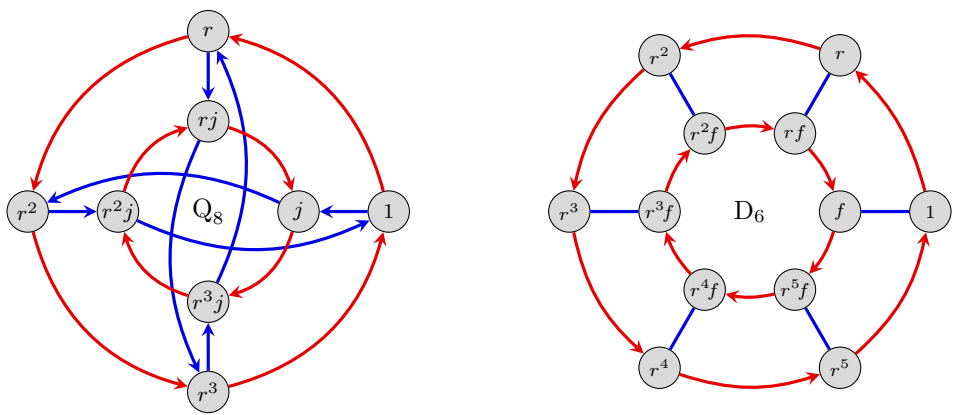

Cayley graphs are basically a visualization of a group presentation, and are useful for quickly multiplying elements. Because they depend on a presentation, different choices of generators can lead to very different looking graphs. For example, the Cayley graph of $\Dih_6=\<s,t\>$, where $s$ and $t$ are two reflections (e.g., $f$ and $rf$), would look like a $12$-gon, because $s^2=t^2=(st)^6=1$. Our second visualization of a group is its \emph{cycle graph}. Like Cayley graphs, the vertex set of a cycle graph consists of the elements in the group. Each one generates a cyclic subgroup (an ``orbit''), and the cycle graph depicts the maximal orbits, and how they intersect. Unlike Cayley graphs, these do not provide information about how to multiply, but they are independent of the generating set, and so all edges are undirected. Cycle graphs for $\Dih_6$ and $\Quat_8$ are shown in Figure~\ref{fig:D4-Q8-cycle}. 

\tikzstyle{v-small} = [circle, draw, fill=lightgrey,inner sep=0pt,
    minimum size=3.5mm]
  \tikzstyle{v-r} = [circle, draw, fill=lred,inner sep=0pt, minimum size=5mm]
  \tikzstyle{v-g} = [circle, draw, fill=lightgrey,inner sep=0pt, minimum size=5mm]
  \tikzstyle{v-y} = [circle, draw, fill=lyellow,inner sep=0pt, minimum size=5mm]
\tikzstyle{v-p} = [circle, draw, fill=conj-purple,inner sep=0pt, minimum size=5mm]
   \tikzstyle{v-b} = [circle, draw, fill=lblue,inner sep=0pt,
    minimum size=5mm]
\tikzstyle{edge} = [draw]
\tikzstyle{cy2} = [draw,thick]
\tikzstyle{cy2-r} = [draw,thick]
\tikzstyle{cy2-rs} = [draw,Blue,thick]
\tikzstyle{cy2-r2s} = [draw,amethyst,thick]
\tikzstyle{cy2-r3s} = [draw,darkgreen,thick]
  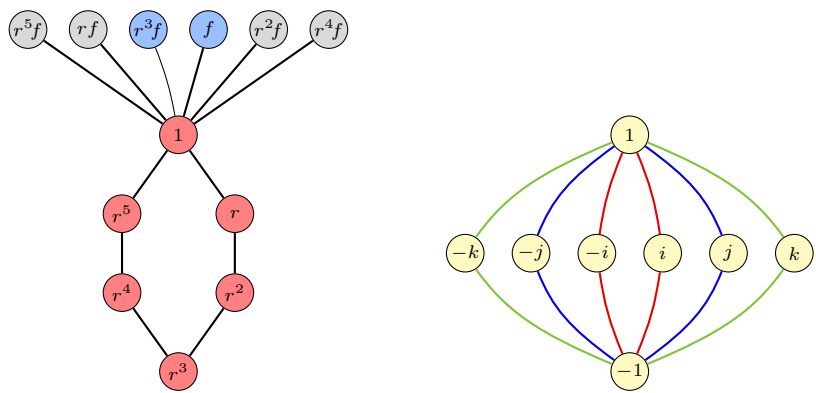
\begin{figure}[!ht]
  \begin{tikzpicture}[scale=1,yscale=.7]
    \tikzstyle{every node}=[font=\scriptsize]
    \begin{scope}[shift={(0,0)}]
      \node (1) at (0,3) [v-r] {$1$};
      \node (r) at (.75,1.5) [v-r] {$r$};
      \node (r2) at (.75,0) [v-r] {$r^2$};
      \node (r3) at (0,-1.5) [v-r] {$r^3$};
      \node (r4) at (-.75,0) [v-r] {$r^4$};
      \node (r5) at (-.75,1.5) [v-r] {$r^5$};
      \node (r5s) at (-2,5) [v-g] {$r^5\!f$};
      \node (rs) at (-1.2,5) [v-g] {$rf$};
      \node (r3s) at (-.4,5) [v-b] {$r^3\!f$};
      \node (s) at (.4,5) [v-b] {$f$};
      \node (r2s) at (1.2,5) [v-g] {$r^2\!f$};
      \node (r4s) at (2,5) [v-g] {$r^4\!f$};
      \draw [cy2] (1) to (r); \draw [cy2] (r) to (r2);
      \draw [cy2] (r2) to (r3); \draw [cy2] (r3) to (r4);
      \draw [cy2] (r4) to (r5); \draw [cy2] (r5) to (1);
      \draw [cy2] (1) to (s); \draw [cy2] (1) to (r2s);
      \draw [cy2] (1) to (r4s); 
      \draw [cy2] (1) to (rs); \draw [edge, bend right=4] (1) to (r3s);
      \draw [cy2] (1) to (r5s); 
    \end{scope}
    \begin{scope}[shift={(6,0)},xscale=1.75]
    \tikzstyle{v-yel} = [circle, draw, fill=lyellow,inner sep=0pt,
    minimum size=5mm]
    \tikzstyle{every node}=[font=\scriptsize]
      \node (1) at (0,3) [v-yel] {$1$};
      \node (-1) at (0,-1.5) [v-yel] {$-1$};
      \node (-k) at (-1.25,.75) [v-yel] {$-k$};
      \node (-j) at (-.75,.75) [v-yel] {$-j$};
      \node (-i) at (-.25,.75) [v-yel] {$-i$};
      \node (i) at (.25,.75) [v-yel] {$i$};
      \node (j) at (.75,.75) [v-yel] {$j$};
      \node (k) at (1.25,.75) [v-yel] {$k$};

      \draw [cy2,Green] (1) to [bend right=14] (-k);
      \draw [cy2,Blue] (1) to [bend right=10] (-j);
      \draw [cy2,Red] (1) to [bend right=3] (-i);
      \draw [cy2,Red] (1) to [bend left=3] (i);
      \draw [cy2,Blue] (1) to [bend left=10] (j);
      \draw [cy2,Green] (1) to [bend left=14] (k);
      \draw [cy2,Green] (-1) to [bend right=14] (k);
      \draw [cy2,Blue] (-1) to [bend right=10] (j);
      \draw [cy2,Red] (-1) to [bend right=3] (i);
      \draw [cy2,Red] (-1) to [bend left=3] (-i);
      \draw [cy2,Blue] (-1) to [bend left=10] (-j);
      \draw [cy2,Green] (-1) to [bend left=14] (-k);
    \end{scope}
    \end{tikzpicture}
    \caption{The cycle graphs of the dihedral group $\Dih_6$ (left) and the quaternion group $\Quat_8$ (right).}\label{fig:D4-Q8-cycle}
    \end{figure}

  Throughout this article, one should continue to think of the generator $r$ as a counterclockwise rotation by $2\pi/n$ radians. At times, we will write this as a primitive $n^{\rm th}$ root of unity $\zeta=\zeta_n=e^{2\pi i/n}$, and in other settings, it is convenient to write it as the canonical $2\times 2$ reflection matrix, as in Equation~\eqref{eq:Dn-representation}. Most of our groups will be generated by $r$, and some other element $s$, which takes the place of $f$ in the dihedral group. 
  
  \section{Dicyclic groups}  
  
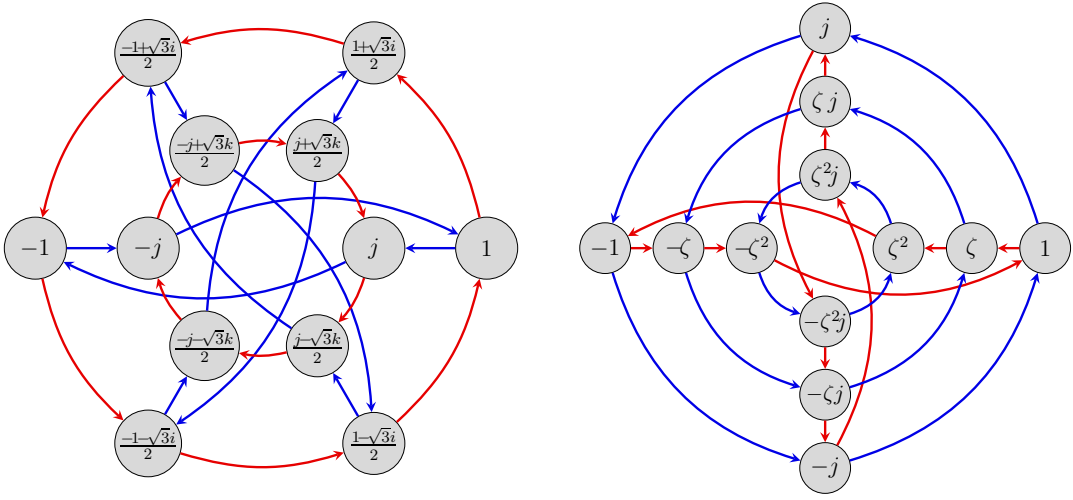
\begin{figure}[!ht]
  \tikzstyle{R-out} = [draw, very thick, Red,-stealth,bend right=17]
  \tikzstyle{R-in} = [draw, very thick, Red,-stealth,bend left=12]
  \tikzstyle{B} = [draw, very thick, Blue,-stealth,bend left=26]
  \scalebox{.75}{
\begin{tikzpicture}[scale=2,baseline=1.3ex,auto]
  \tikzstyle{v-big} = [circle,fill=lightgrey,draw,inner sep=0pt, minimum size=11mm]
    \tikzstyle{v} = [circle,fill=lightgrey,draw,inner sep=0pt, minimum size=9mm]
\tikzstyle{every node}=[font=\Large]
    \begin{scope}[shift={(0,0)},scale=1]
        \node (1) at (0:2.0) [v-big] {$1$};
        \node (r) at (60:2.0) [v-big] {$\frac{1\!+\!\sqrt{3}i}{2}$};
        \node (r2) at (120:2.0) [v-big] {\large $\frac{-\!1\!+\!\sqrt{3}i}{2}$};
        \node (r3) at (180:2.0) [v-big] {$-1$};
        \node (r4) at (240:2.0) [v-big] {\large $\frac{-\!1\!-\!\sqrt{3}i}{2}$};
        \node (r5) at (300:2.0) [v-big] {$\frac{1\!-\!\sqrt{3}i}{2}$};
        \node (s) at (0:1) [v-big] {$j$};
        \node (sr5) at (60:1) [v-big] {$\frac{j\!+\!\sqrt{3}k}{2}$};
        \node (sr4) at (120:1) [v-big] {\large $\frac{-\!j\!+\!\sqrt{3}k}{2}$};
        \node (sr3) at (180:1) [v-big] {$-j$};
        \node (sr2) at (240:1) [v-big] {\large $\frac{-\!j\!-\!\sqrt{3}k}{2}$};
        \node (sr) at (300:1) [v-big] {$\frac{j\!-\!\sqrt{3}k}{2}$};
        \path[b] (1) to (s);
        \path[B] (s) to (r3);
        \path[b] (r3) to (sr3);
        \path[B] (sr3) to (1);
        \path[b] (r) to (sr5);
        \path[B] (sr5) to (r4);
        \path[b] (r4) to (sr2);
        \path[B] (sr2) to (r);
        \path[b] (r2) to (sr4);
        \path[B] (sr4) to (r5);
        \path[b] (r5) to (sr);
        \path[B] (sr) to (r2);
        \path[R-out] (1) to (r);
        \path[R-out] (r) to (r2);
        \path[R-out] (r2) to (r3);
        \path[R-out] (r3) to (r4);
        \path[R-out] (r4) to (r5);
        \path[R-out] (r5) to (1);
        \path[R-in] (sr5) to (s);
        \path[R-in] (sr4) to (sr5);
        \path[R-in] (sr3) to (sr4);
        \path[R-in] (sr2) to (sr3);
        \path[R-in] (sr) to (sr2);
        \path[R-in] (s) to (sr);
    \end{scope}
    \begin{scope}[shift={(5,0)},scale=.65]
      \tikzstyle{every node}=[font=\Large]
      \node[v] (1) at (0:3) {$1$};
      \node[v] (r) at (0:2) {$\zeta$};
      \node[v] (r2) at (0:1) {$\zeta^2$};
      \node[v] (s) at (90:3) {$j$};
      \node[v] (rs) at (90:2) {$\zeta\,j$};
      \node[v] (r2s) at (90:1) {$\zeta^2\!j$};
      \node[v] (s2) at (180:3) {$-1$};
      \node[v] (rs2) at (180:2) {$-\!\zeta$};
      \node[v] (r2s2) at (180:1) {$-\!\zeta^2$};
      \node[v] (s3) at (-90:3) {$-j$};
      \node[v] (rs3) at (-90:2) {\large $-\zeta j$};
      \node[v] (r2s3) at (-90:1) {\large $-\zeta^2\!j$};
      \draw[r] (1) to (r);
      \draw[r] (r) to (r2);
      \draw[r] (s2) to (rs2);
      \draw[r] (rs2) to (r2s2);
      \draw[r] (r2s) to (rs);
      \draw[r] (rs) to (s);
      \draw[r] (r2s3) to (rs3);
      \draw[r] (rs3) to (s3);
      \draw[r,bend right=28] (r2) to (s2); 
      \draw[r,bend right=28] (r2s2) to (1); 
      \draw[r,bend right=28] (s) to (r2s3); 
      \draw[r,bend right=28] (s3) to (r2s); 
      \draw[b,bend right=28] (1) to (s);
      \draw[b,bend right=28] (r) to (rs);
      \draw[b,bend right=28] (r2) to (r2s);
      \draw[b,bend right=28] (s) to (s2);
      \draw[b,bend right=28] (rs) to (rs2);
      \draw[b,bend right=28] (r2s) to (r2s2);
      \draw[b,bend right=28] (s2) to (s3);
      \draw[b,bend right=28] (rs2) to (rs3);
      \draw[b,bend right=28] (r2s2) to (r2s3);
      \draw[b,bend right=28] (s3) to (1);
      \draw[b,bend right=28] (rs3) to (r);
      \draw[b,bend right=28] (r2s3) to (r2);
    \end{scope}
\end{tikzpicture}}
\caption{Two ways to lay out a Cayley graph for the dicyclic group $\Dic_6=\<\zeta_6,j\>$.}\label{fig:Dic6-Cayley}
\end{figure}

  Our first family of lesser-known groups has a simple construction: start with $\Quat_8=\<i,j\>$, and replace $i=\zeta_4=e^{2\pi i/4}$ with a larger (even) root of unity, $\zeta_n=\cos(2\pi/n)+i\sin(2\pi/n)$. The multiplication rules, like $ij=k$ and $jk=i$, remain unchanged. This defines the \emph{dicyclic group},
  \[
  \Dic_n=\big\<\zeta_n,j\big\>\cong\left\<\underbrace{\Alert{\begin{bmatrix}\zeta_n & 0 \\ 0 & \overline{\zeta}_n\end{bmatrix}}}_{R=R_n},\;
  \underbrace{\Balert{\begin{bmatrix}0 & -1 \\ 1 & 0\end{bmatrix}}}_S\right\>.
  \]
  Figure~\ref{fig:Dic6-Cayley} shows two Cayley graphs for $\Dic_6$. On the left, the roots of unity are arranged around the outer circle, and their products with $j$ are on the inner circle. On the right, powers of $\zeta$ increase as one moves inward, and the powers of $j$ increase as one traverses around in the positive (counterclockwise) direction. Note that $\zeta^3=-1$ for a sixth root of unity, and in general, $\zeta^{n/2}=-1$. Sometimes, we prefer to write the roots of unity as $1,\zeta,\dots,\zeta^{n-1}$, and other times, $\pm 1,\pm\zeta,\dots,\pm \zeta^{n/2-1}$ is more convenient.
 
One standard presentation of the dicyclic group is motivated by the ``unit circle'' Cayley graph on the left in Figure~\ref{fig:Dic6-Cayley}. Letting $r=\zeta_n$ and $s=j$, we get
\[
 \Dic_n=\big\<r,s\mid r^n=s^4=1,\;rsr=s,\;r^{n/2}=s^2\big\>.
\]
Notice how the special case of $n=4$ yields the quaternion group $\Quat_8$ from Eq.~\eqref{eq:Q8-representation}. With a Cayley graph, certain properties and features become immediately apparent. For example, we can see that $\Dic_6$ is generated by $\zeta^2$ and $j$. This leads to an alternative presentation, $\Dic_6=\<a,b\mid a^4=b^3,\;bab=a\>$. 

Our notation for the dicyclic groups is nonstandard. Most sources refer to the group $\Dic_6$ as $\Dic_3$, for ``\emph{the third dicyclic group}.'' This convention explains the ``di-'' prefix, because it means that $\Dic_n$ is an extension of $\Dih_n$ by $\Cyc_{2n}$, via $1\to \Cyc_{2n}\to\Dic_n \to \Dih_n\to 1$. Other sources use $\Dic_{12}$, for ``\emph{the dicyclic group with $12$ elements}.'' We strongly feel that the subscript should reflect the specific roots of unity used, and the remainder of this paper should strengthen that argument.

Using our convention that $\Dic_n=\<\zeta_n,j\>=\<r,s\>$ is a group of order $2n$, it is an extension of the dihedral group $\Dih_{n/2}$ by $\Cyc_n$, via a short exact sequence $1\to \Cyc_n\to\Dic_n \to \Dih_{n/2}\to 1$. The isomorphism $\Dic_6/\<-1\>\cong \Dih_3$ should be apparent from the subgroup lattice and Cayley graph shown in Figure~\ref{fig:Dic6-quotient}. By the correspondence theorem, the lattice of the quotient $\Dic_6/\<r^3\>$ (recall that $r^3=-1$) has the same structure as the intermediate sublattice between $\Dic_6$ and $\<r^3\>$, which is highlighted in Figure~\ref{fig:Dic6-quotient}. These types of observations are compelling evidence for the utility of our last visual tool, \emph{subgroup lattices}. 

\tikzstyle{v-big} = [circle,fill=lightgrey,draw,inner sep=0pt, minimum size=6mm]
\begin{figure}[!ht]
  \begin{tikzpicture}[scale=1.15,box/.style={anchor=south}]
    \begin{scope}[shift={(0,0)},scale=.6]  
     \tikzstyle{every node}=[font=\scriptsize]
     \tikzstyle{v-r} = [circle, draw, fill=lred,inner sep=0pt,
    minimum size=5mm]
   \tikzstyle{v-b} = [circle, draw, fill=lblue,inner sep=0pt,
    minimum size=5mm]
    \tikzstyle{v-p} = [circle, draw, fill=lightpurple,inner sep=0pt, 
    minimum size=5mm]
         \tikzstyle{v-small} = [circle, draw, fill=lightgrey,inner sep=0pt,
    minimum size=5mm]
      \tikzstyle{cy2} = [draw]
  \tikzstyle{cy2-r} = [draw]
  \tikzstyle{cy2-s} = [draw,Blue]
  \tikzstyle{cy2-rs} = [draw,Green]
  \tikzstyle{cy2-r2s} = [draw,Purple]
  \tikzstyle{cy2-r3s} = [draw,orange] 
      \node (1) at (0,3) [v-r] {$1$};
      \node (s2) at (0,-3) [v-r] {$s^2$};
      \node (r) at (.5,1) [v-r] {$r$};
      \node (r2) at (.5,-1) [v-r] {$r^2$};
      \node (rs2) at (-.5,-1) [v-r] {$r^4$};
      \node (r2s2) at (-.5,1) [v-r] {$r^5$};
      \node (s3) at (-1.5,0) [v-b] {$s^3$};
      \node (s) at (1.5,0) [v-b] {$s$};
      \node (rs3) at (-2.5,0) [v-small] {$rs^3$};
      \node (rs) at (2.5,0) [v-small] {$rs$};
      \node (r2s3) at (-3.5,0) [v-small] {\tiny $r^2\!\!s^3$};
      \node (r2s) at (3.5,0) [v-small] {$r^2\!s$};
      \draw [cy2] (1) to (r);
      \draw [cy2] (r) to (r2);
      \draw [cy2] (r2) to (s2);
      \draw [cy2] (s2) to (rs2);
      \draw [cy2] (rs2) to (r2s2);
      \draw [cy2] (r2s2) to (1);
      \draw [cy2-s,bend left=10] (1) to (s);
      \draw [cy2-s,bend right=10] (s2) to (s);
      \draw [cy2-rs,bend left=5] (1) to (rs);
      \draw [cy2-rs,bend right=5] (s2) to (rs);
      \draw [cy2-r2s] (1) to (r2s); \draw [cy2-r2s] (s2) to (r2s);
      \draw [cy2-s, bend right=10] (1) to (s3);
      \draw [cy2-s,bend left=10] (s2) to (s3);
      \draw [cy2-rs,bend right=5] (1) to (rs3);
      \draw [cy2-rs,bend left=5] (s2) to (rs3);
      \draw [cy2-r2s] (1) to (r2s3); \draw [cy2-r2s] (s2) to (r2s3);
    \end{scope}
    \begin{scope}[shift={(4.25,-2.1)},scale=.72,shorten >= -2pt, shorten <= -2pt,xscale=.8,yscale=.95]
      \tikzstyle{every node}=[font=\normalsize]
      \node(G) at (0,6) {$\Dic_6$};
      \node(b) at (-1.5,4.5) {$\<r\>$};
      \node(aa) at (0,1.5) {$\<r^3\>$};
      \node[faded] (bb) at (-2,2.25) {$\<r^2\>$};
      \node(a) at (.5,24/7) {$\<s\>$};
      \node(ab) at (1.75,24/7) { $\<rs\>$};
      \node(ba) at (3,24/7) {$\<r^2\!s\>$};
      \node[faded] (1) at (0,0) {$\left\<1\right\>$};
      \draw [faded] (1) to (aa);
      \draw [faded] (1) to (bb);
      \draw (b) to (aa);
      \draw [faded] (b) to (bb);
      \draw (aa) to (a);
      \draw (aa) to (ab);
      \draw (aa) to (ba);
      \draw (G) to (a);
      \draw (G) to (ab);
      \draw (G) to (ba);
      \draw (G) to (b);   
    \end{scope}
    \begin{scope}[shift={(8.25,0)},scale=1,scale=.8]
    \tikzstyle{every node}=[font=\footnotesize]
      \node (1) at (0:2.2) [v-big] {$\pm 1$};
      \node (r) at (120:2.2) [v-big] {$\pm\zeta$};
      \node (r2) at (240:2.2) [v-big] {$\pm\zeta^2$};
      \node (f) at (0:1) [v-big] {$\pm j$};
      \node (r2f) at (240:1) [v-big] {\scriptsize $\pm\zeta^2\!j$};
      \node (rf) at (120:1) [v-big] {$\pm\zeta j$};
      \draw [r] (1) to [bend right] (r);
      \draw [r] (r) to [bend right] (r2);
      \draw [r] (r2) to [bend right] (1);
      \draw [r] (f) to [bend left] (r2f);
      \draw [r] (r2f) to [bend left] (rf);
      \draw [r] (rf) to [bend left] (f);
      \draw [bb] (1) to (f);
      \draw [bb] (r) to (rf);
      \draw [bb] (r2) to (r2f);
    \end{scope}
  \end{tikzpicture}
  \caption{Left: the cycle graph of the dicyclic group $\Dic_6$, which generalizes the cycle graphs of both $\Dih_6$ and $\Quat_8$ in Figure~\ref{fig:D4-Q8-cycle}. Middle: The structure of the subgroup lattice of $\Dih_3$ appears at the top of the lattice of $\Dic_6$ because it is a quotient. Right: The Cayley graph showing $\Dic_6\!/\<r^3\>$; recall $r^3=-1$.}\label{fig:Dic6-quotient}
\end{figure}
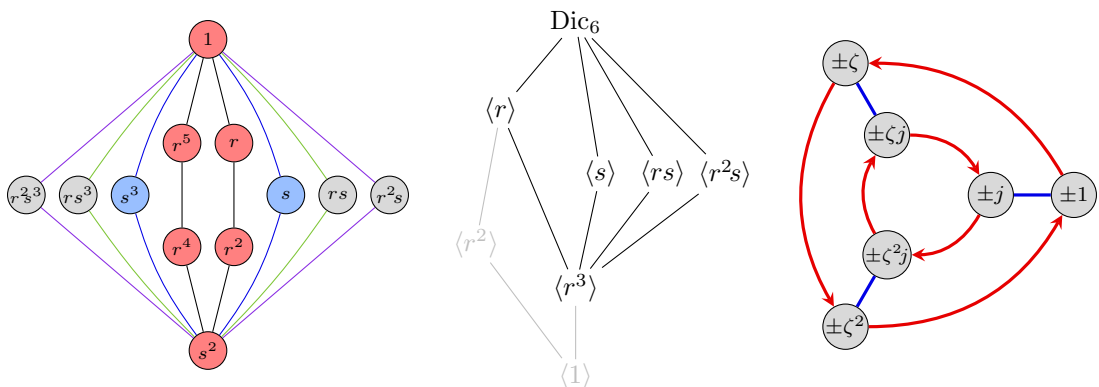

The cycle graph of $\Dic_6$ also appears in Figure~\ref{fig:Dic6-quotient}, and it is illuminating to compare it to those of $\Dih_6$ and $\Quat_8$, the groups that we put together to build it. Notice how in $\Dih_6$, the six elements not in the subgroup $\<r\>$ have order $2$ because they are reflections. In comparison, in $\Dic_6$, they have order $4$, because they are multiples of $j$. The generalization of $\Dic_6$ from $\Quat_8$ is also clear from the cycle graphs -- simply replace the $4^{\rm th}$ roots of unity with the $6^{\rm th}$ roots.

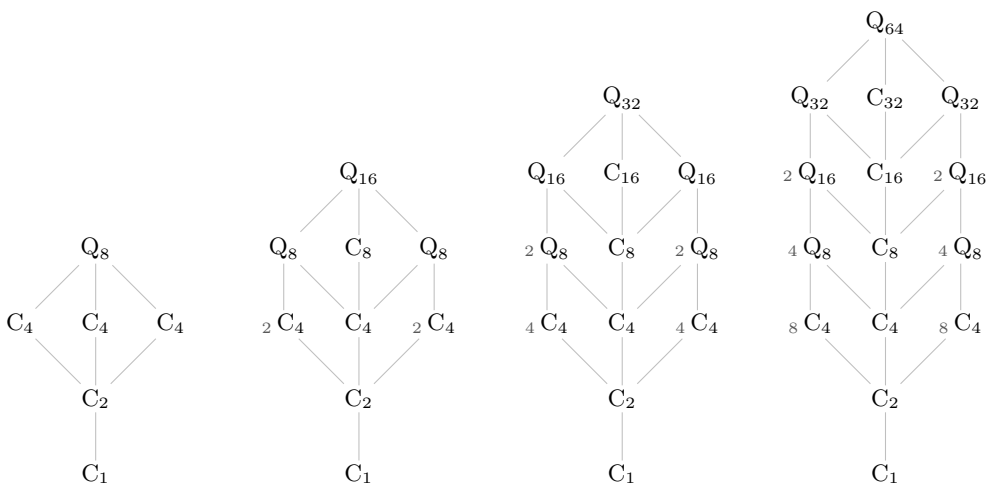
\begin{figure}[!ht]
  \begin{tikzpicture}[scale=1,shorten >= -2pt, shorten <= -2pt]
    \tikzstyle{every node}=[font=\small]
   \begin{scope}[shift={(0,0)}]
     \node (Q8) at (0,3) {$\Quat_8$};
     \node (C4) at (0,2) {$\Cyc_4$};
     \node (C2) at (0,1) {$\Cyc_2$};
     \node (1) at (0,0) {$\Cyc_1$};
     \node (C4-l) at (-1,2) {$\Cyc_4$};
     \node (C4-r) at (1,2) {$\Cyc_4$};
     \draw[f] (Q8) to (C4);
     \draw[f] (C4) to (C2);
     \draw[f] (C2) to (1);
     \draw[f] (Q8) to (C4-l); \draw[f] (Q8) to (C4-r); 
     \draw[f] (C4-l) to (C2); \draw[f] (C4-r) to (C2);
  \end{scope}
\begin{scope}[shift={(3.5,0)}]
     \node (Q16) at (0,4) {$\Quat_{16}$};
     \node (C8) at (0,3) {$\Cyc_8$};
     \node (C4) at (0,2) {$\Cyc_4$};
     \node (C2) at (0,1) {$\Cyc_2$};
     \node (1) at (0,0) {$\Cyc_1$};
     \node (Q8-l) at (-1,3) {$\Quat_8$};
     \node (C4-l) at (-1,2) {${}_{\color{midgrey} 2}\Cyc_4$};
     \node (Q8-r) at (1,3) {$\Quat_8$};
     \node (C4-r) at (1,2) {${}_{\color{midgrey} 2}\Cyc_4$};
     \draw[f] (Q16) to (C8);
     \draw[f] (C8) to (C4);
     \draw[f] (C4) to (C2);
     \draw[f] (C2) to (1);
     \draw[f] (Q16) to (Q8-l); \draw[f] (Q16) to (Q8-r); 
     \draw[f] (Q8-l) to (C4-l); \draw[f] (Q8-r) to (C4-r); 
     \draw[f] (Q8-l) to (C4); \draw[f] (Q8-r) to (C4);
     \draw[f] (C4-l) to (C2); \draw[f] (C4-r) to (C2);
  \end{scope}
 \begin{scope}[shift={(7,0)}]
     \node (Q32) at (0,5) {$\Quat_{32}$};
     \node (C16) at (0,4) {$\Cyc_{16}$};
     \node (C8) at (0,3) {$\Cyc_8$};
     \node (C4) at (0,2) {$\Cyc_4$};
     \node (C2) at (0,1) {$\Cyc_2$};
     \node (1) at (0,0) {$\Cyc_1$};
     \node (Q16-l) at (-1,4) {$\Quat_{16}$};
     \node (Q8-l) at (-1,3) {${}_{\color{midgrey} 2}\Quat_8$};
     \node (C4-l) at (-1,2) {${}_{\color{midgrey} 4}\Cyc_4$};
     \node (Q16-r) at (1,4) {$\Quat_{16}$};
     \node (Q8-r) at (1,3) {${}_{\color{midgrey} 2}\Quat_8$};
     \node (C4-r) at (1,2) {${}_{\color{midgrey} 4}\Cyc_4$};
     \draw[f] (Q32) to (C16);
     \draw[f] (C16) to (C8);
     \draw[f] (C8) to (C4);
     \draw[f] (C4) to (C2);
     \draw[f] (C2) to (1);
     \draw[f] (Q32) to (Q16-l); \draw[f] (Q32) to (Q16-r); 
     \draw[f] (Q16-l) to (Q8-l); \draw[f] (Q16-r) to (Q8-r); 
     \draw[f] (Q8-l) to (C4-l); \draw[f] (Q8-r) to (C4-r); 
     \draw[f] (Q16-l) to (C8); \draw[f] (Q16-r) to (C8);
     \draw[f] (Q8-l) to (C4); \draw[f] (Q8-r) to (C4);
     \draw[f] (C4-l) to (C2); \draw[f] (C4-r) to (C2);
  \end{scope}
  \begin{scope}[shift={(10.5,0)}]
     \node (Q64) at (0,6) {$\Quat_{64}$};
     \node (C32) at (0,5) {$\Cyc_{32}$};
     \node (C16) at (0,4) {$\Cyc_{16}$};
     \node (C8) at (0,3) {$\Cyc_8$};
     \node (C4) at (0,2) {$\Cyc_4$};
     \node (C2) at (0,1) {$\Cyc_2$};
     \node (1) at (0,0) {$\Cyc_1$};
     \node (Q32-l) at (-1,5) {$\Quat_{32}$};
     \node (Q16-l) at (-1,4) {${}_{\color{midgrey} 2}\Quat_{16}$};
     \node (Q8-l) at (-1,3) {${}_{\color{midgrey} 4}\Quat_8$};
     \node (C4-l) at (-1,2) {${}_{\color{midgrey} 8}\Cyc_4$};
     \node (Q32-r) at (1,5) {$\Quat_{32}$};
     \node (Q16-r) at (1,4) {${}_{\color{midgrey} 2}\Quat_{16}$};
     \node (Q8-r) at (1,3) {${}_{\color{midgrey} 4}\Quat_8$};
     \node (C4-r) at (1,2) {${}_{\color{midgrey} 8}\Cyc_4$};
     \draw[f] (Q64) to (C32);
     \draw[f] (C32) to (C16);
     \draw[f] (C16) to (C8);
     \draw[f] (C8) to (C4);
     \draw[f] (C4) to (C2);
     \draw[f] (C2) to (1);
     \draw[f] (Q64) to (Q32-l); \draw[f] (Q64) to (Q32-r);
     \draw[f] (Q32-l) to (Q16-l); \draw[f] (Q32-r) to (Q16-r); 
     \draw[f] (Q16-l) to (Q8-l); \draw[f] (Q16-r) to (Q8-r);
     \draw[f] (Q8-l) to (C4-l); \draw[f] (Q8-r) to (C4-r); 
     \draw[f] (Q32-l) to (C16); \draw[f] (Q32-r) to (C16);
     \draw[f] (Q16-l) to (C8); \draw[f] (Q16-r) to (C8);
     \draw[f] (Q8-l) to (C4); \draw[f] (Q8-r) to (C4);
     \draw[f] (C4-l) to (C2); \draw[f] (C4-r) to (C2);
  \end{scope}
\end{tikzpicture}
\caption{The ``reduced'' subgroup lattices for the generalized quaternion groups $\Quat_{2^n}=\Dic_{2^{n-1}}$ follow a predictable pattern. Each group shown represents a conjugacy class, and the left-subscripts denote its size. }\label{fig:Qn-conjugacy}
\end{figure}

In this article, we are keeping the name ``dicyclic,'' despite deviating from the standard subscript that inspired the ``di-'' prefix in the first place. Some books call these groups the \emph{generalized quaternion} groups, and use the subscript to denote their order. For example, the group above would be $\Quat_{12}$. We will take a more common approach, and only use the term ``generalized quaternion'' for dicyclic groups whose order is a power of $2$. Our reasoning is because dicyclic groups of order $2^m$ more closely generalize the structure of the classic quaternion group $\Quat_8$ than the other dicyclic groups do. One notable property of $\Quat_8$ is that any two nontrivial cyclic subgroups contain $-1$. This means that nontrivial subgroups intersect nontrivially, and so $\Quat_8$ does not break up into a semidirect product of its proper subgroups. When $n=2^m$, every $n^{\rm th}$ root of unity generates a cyclic subgroup that contains $r^{n/2}=-1$. Therefore, any two nontrivial subgroups of $\Dic_n$ also intersect in at least $\<-1\>$, and so these dicyclic groups do not decompose as semidirect products either. This is one reason why only these groups are typically considered ``generalized quaternion.'' This property, along with the fact that they have a dihedral quotient, is apparent in the subgroup lattices in Figure~\ref{fig:Qn-conjugacy}. The nodes here are actually conjugacy classes of subgroups, with the left subscript denoting the size.\footnote{In general, a subgroup lattice collapsed in this manner need not be an actual lattice.} This convention follows the online group databases GroupNames \cite{GroupNames} and LMFDB \cite{LMFDB}.

The pattern of the generalized quaternion groups from Figure~\ref{fig:Qn-conjugacy} should be clear -- their subgroup lattices look like those of ``\emph{the dihedral group $\Dih_{2^{n-2}}$ on a stick}.'' By the correspondence theorem, $\Quat_{2^n}/\<-1\>$ is isomorphic to $\Dih_{2^{n-2}}$. Compare this to the quotient $\Dic_6/\<-1\>\cong \Dih_3$ shown in   Figure~\ref{fig:Dic6-quotient}. In that group, $\<r^2\>\cong \Cyc_3$ is normal, its product with the non-normal subgroup $\<s\>\cong \Cyc_4$ is $\Dic_6$, and these subgroups intersect trivially. Therefore, $\Dic_6\cong \Cyc_3\rtimes \Cyc_4$. 

Our last comment on the dicyclic groups is that they are only defined for even $n$. Of course, one can construct such a group using an odd root of unity $\zeta_n=\zeta_{2k+1}$, but that ends up generating $\zeta_{2n}$, and hence the dicyclic group $\Dic_{2n}$. Therefore, it is standard to assume that $n$ is even.

\section{Diquaternion groups}
 
In the previous section, we constructed the dicyclic groups by starting with $\Quat_8=\<\zeta_4,j\>$, and replacing $i=\zeta_4$ with a larger root of unity, which can be thought of as a generator (rotation) from $\Dih_n$. Another way we can combine these groups, which is more reflective of the idea of ``\emph{dihedralizing the quaternions},'' is to start with $\Quat_8$ and throw in the other generator of $\Dih_n$ --- a reflection, $f$. Though it may not be clear what this means in terms of generators and relations, it is unambiguous using our matrix representations. We will call this the \emph{diquaternion group}, and write
  \[
  \DQ_8=\big\<i,j,k,f\big\>\cong\left\<
  \underbrace{\begin{bmatrix} i & 0 \\ 0 & -i\end{bmatrix}}_{R=R_4},
  \underbrace{\begin{bmatrix} 0 & -1 \\ 1 & 0\end{bmatrix}}_{S},
  \underbrace{\begin{bmatrix} 0 & -i \\ -i & 0\end{bmatrix}}_{T},
  \underbrace{\Balert{\begin{bmatrix} 0 & 1 \\ 1 & 0\end{bmatrix}}}_{F}\right\>.
  \]
  Multiplying any $2\times 2$ matrix by the ``complex reflection marix'' $F$ swaps its rows or columns, depending on whether it is left- or right-multiplied. Thus, the result of including $F$ is that it doubles the size of $\Quat_8$ --- each matrix has a corresponding ``mirror reflection.'' Of course, the generator $k$ (and matrix $T$) above is unnecessary, but it is helpful to include. In quantum physics, this group is known as the 
  \emph{Pauli group on $1$ qubit}, and is 
  \[
  \DQ_8=\big\<X,Y,Z\big\>=\big\{\pm\!I,\,\pm i I,\,\pm X,\,\pm iX,\,\pm Y,\,\pm iY,\,\pm Z,\,\pm iZ\big\},
  \]
  generated by the \emph{Pauli matrices} from quantum mechanics and
  information theory:
  \[
  X=\begin{bmatrix}0&1\\1&0\end{bmatrix},\qquad
  Y=\begin{bmatrix}0&-i\\ i&0\end{bmatrix},\qquad
  Z=\begin{bmatrix}1&0\\0&-1\end{bmatrix}.
  \]
  It is straightforward to check that our two generating sets for $\DQ_8$ are related by
  \[
  XY=R\quad\text{``$i$''},\qquad\qquad
  XZ=S\quad\text{``$j$''},\qquad\qquad
  YZ=\overline{T}\quad\text{``$-k$''},
  \]
  where $T=RS$. Though the diquaternion group is better motivated by adding the reflection matrix $F$ to $\Quat_8\cong \<R,S\>$, the Cayley graph is arguably more pleasing using the Pauli matrices as generators, which is shown in Figure~\ref{fig:DQ8-Cayley}. Nodes with matrices corresponding to the quaternions are labeled and highlighted.

  \tikzstyle{v-small} = [circle, draw, fill=lightgrey,inner sep=0pt, 
  minimum size=7mm]
\tikzstyle{v-tiny} = [circle, draw, fill=lightgrey,inner sep=0pt, 
  minimum size=1.25mm]
    \tikzstyle{v-yel} = [circle, draw, fill=lyellow,inner sep=0pt,
      minimum size=7mm]
\tikzstyle{rr-bend} = [draw, very thick, Red,bend left=12]
\tikzstyle{bb-bend} = [draw, very thick, Blue,bend left=12]
\begin{figure}[!ht]
\begin{tikzpicture}[scale=1.7,auto]
  \tikzstyle{every node}=[font=\tiny]
  \tikzstyle{rr-bend} = [draw, very thick, Red,bend left=12]
  \tikzstyle{bb-bend} = [draw, very thick, Blue,bend left=12]
  \tikzstyle{rr} = [draw, very thick, Red]
  \tikzstyle{bb} = [draw, very thick, Blue]
  \tikzstyle{gg} = [draw, ultra thick, Green]
  \begin{scope}[shift={(0,0)}]
    \node at (-4,1.5)
          {\normalsize\Balert{$X=F=\begin{bmatrix}0&1\\1&0\end{bmatrix}$}};
     \node at (-4,.5)
                {\normalsize\Alert{$Y:=\begin{bmatrix}0&-i\\i&0\end{bmatrix}$}};
    \node at (-4,-.5)
                      {\normalsize\Galert{$Z=\begin{bmatrix}1&0\\0&-1\end{bmatrix}$}};
    \node at (-4,-1.5) {\normalsize $\Balert{X}\Alert{Y}=R$,\quad $\Balert{X}\Galert{Z}=S$,\quad $\Alert{Y}\Galert{Z}=\overline{T}$};
    \setlength{\arraycolsep}{1.5pt}
    \renewcommand{\arraystretch}{.7}
      \node (s) at (90:1.15) [v-small] {$\begin{bmatrix}i&0\\0&i\end{bmatrix}$};
      \node (rs) at (45:1.15) [v-yel] {$\begin{bmatrix}0&-\!1\\1&0\end{bmatrix}$};
      \node (r2s) at (0:1.15) [v-small] {$\begin{bmatrix}1&0\\0&-\!1\end{bmatrix}$};
      \node (r3s) at (-45:1.15) [v-yel] {$\begin{bmatrix}0&i\\i&0\end{bmatrix}$};
      \node (r4s) at (-90:1.15) [v-small] {$\begin{bmatrix}-\!i&0\\0&\!-\!i\end{bmatrix}$};
      \node (r5s) at (-135:1.15) [v-yel] {$\begin{bmatrix}0&1\\-\!1&0\end{bmatrix}$};
      \node (r6s) at (180:1.15) [v-small] {$\begin{bmatrix}-\!1&0\\0&1\end{bmatrix}$};
      \node (r7s) at (135:1.15) [v-yel] {$\begin{bmatrix}0&-\!i\\\!\!-\!i&0\end{bmatrix}$};
      \node (1) at (90:2.1) [v-yel] {$\begin{bmatrix}i&0\\0&-\!i\end{bmatrix}$};
      \node (r) at (45:2.1) [v-small] {$\begin{bmatrix}0&1\\1&0\end{bmatrix}$};
      \node (r2) at (0:2.1) [v-yel] {$\begin{bmatrix}1&0\\0&1\end{bmatrix}$};
      \node (r3) at (-45:2.1) [v-small] {$\begin{bmatrix}0&-\!i\\i&0\end{bmatrix}$};
      \node (r4) at (-90:2.1) [v-yel] {$\begin{bmatrix}-\!i&0\\0&i\end{bmatrix}$};
      \node (r5) at (-135:2.1) [v-small] {\setlength{\arraycolsep}{1.1pt}$\begin{bmatrix}0&\!\!-\!1\\-\!1&0\end{bmatrix}$};
      \node (r6) at (180:2.1) [v-yel] {\setlength{\arraycolsep}{1.1pt}$\begin{bmatrix}-\!1&0\\0&\!\!-\!1\end{bmatrix}$};
      \node (r7) at (135:2.1) [v-small] {$\begin{bmatrix}0&i\\-\!i&0\end{bmatrix}$};
      \node at (9:2.38) {\large\Palert{$\mathbf{1}$}};
      \node at (90+8:2.4) {\large\Palert{$\mathbf{i}$}};
      \node at (90+28:1.4) {\large\Palert{$\mathbf{k}$}};
      \node at (90-28:1.4) {\large\Palert{$\mathbf{j}$}};
      \node at (180-12:2.38) {\large\Palert{$\mathbf{-1}$}};
      \node at (270-11:2.4) {\large\Palert{$\mathbf{-i}$}};
      \node at (270-30:1.4) {\large\Palert{$\mathbf{-j}$}};
      \node at (270+30:1.4) {\large\Palert{$\mathbf{-k}$}};
      \draw [rr-bend] (1) to (r);
      \draw [bb-bend] (r) to (r2);
      \draw [rr-bend] (r2) to (r3);
      \draw [bb-bend] (r3) to (r4);
      \draw [rr-bend] (r4) to (r5);
      \draw [bb-bend] (r5) to (r6);
      \draw [rr-bend] (r6) to (r7);
      \draw [bb-bend] (r7) to (1);
      \draw [bb] (s) to (r3s);
      \draw [rr] (r3s) to (r6s);
      \draw [bb] (r6s) to (rs);
      \draw [rr] (rs) to (r4s);
      \draw [bb] (r4s) to (r7s);
      \draw [rr] (r7s) to (r2s);
      \draw [bb] (r2s) to (r5s);
      \draw [rr] (r5s) to (s);
      \draw [gg] (1) to (s); \draw [gg] (r) to (rs);
      \draw [gg] (r2) to (r2s); \draw [gg] (r3) to (r3s);
      \draw [gg] (r4) to (r4s); \draw [gg] (r5) to (r5s);
      \draw [gg] (r6) to (r6s); \draw [gg] (r7) to (r7s);
      \node at (0,0) {\large $\DQ_8$};
    \end{scope}
  \end{tikzpicture}
  \caption{The diquaternion group is constructed by adding a reflection matrix $F$ from the dihedral group $\Dih_n$ to the quaternion group $\Quat_8$. Here it is generated by the Pauli matrices from quantum physics. The matrices from the quaternion group are highlighted.}\label{fig:DQ8-Cayley}
  \end{figure}
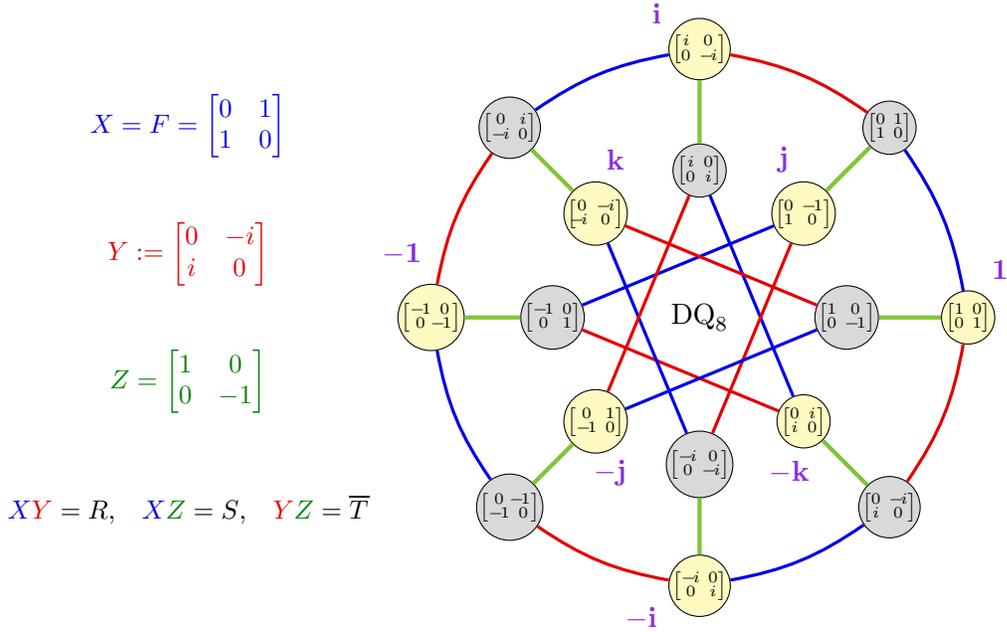
  
  The diquaternion group does not have a simple or intuitive presentation. An example includes the standard but uninspiring
  \begin{equation}\label{eq:DQ8-presentation}
\DQ_8=\big\<a,b,c\mid a^4=c^2=1,\;a^2=b^2,\,ab=ba,\,ac=ca,\,cba=a^2b\big\>.
  \end{equation}
  Even though a Cayley graph does not encode any more information than a presentation, it presents it in a very intuitive format that brings certain structural features to light. For example, from the Cayley graph of $\DQ_8$ in Figure~\ref{fig:DQ8-Cayley}, we can see that there is an index-$2$ subgroup isomorphic to $\Dih_4$, consisting of the nodes along the outer ring. That is, $\<X,Y\>\cong \Dih_4$, where the generators represent noncommuting reflections. We can also see that the nodes along the central axes comprise an abelian subgroup. Specifically, $\<XY,Z\>\cong \Cyc_4\times \Cyc_2$. For both of these (normal) index-$2$ subgroups, there is an element in $\DQ_8$ of order $2$ not in them, and hence a cyclic subgroup $\Cyc_2$ that intersects them trivially. The subgroup $\Quat_8\cong\<XY,XZ\>$ also has such an element. Therefore, $\DQ_8$ decomposes as a semidirect product of each one of these subgroups with $\Cyc_2$. In a classroom, this is a fantastic example of a small group that can be written as a semidirect product in three distinct ways. This example is immediate from inspecting the Cayley graph, but would require a lot of uninspiring work without it, using only the presentation from Equation~\eqref{eq:DQ8-presentation}.
  
  For many groups, the cycle graph highlights certain structural features that Cayley graphs hide, and the dicyclic and diquaternion groups are no exceptions. Their cycle graphs are shown in Figure~\ref{fig:Q16-DQ8-cycle-graphs}, and should be compared to those of the groups they generalize: $\Dih_6$ and $\Quat_8$, back in Figure~\ref{fig:D4-Q8-cycle}.

  \tikzstyle{v-small} = [circle, draw, fill=lightgrey,inner sep=0pt,
    minimum size=3.5mm]
    \tikzstyle{v-r} = [circle, draw, fill=lred,inner sep=0pt,
    minimum size=3.5mm]
   \tikzstyle{v-b} = [circle, draw, fill=lblue,inner sep=0pt,
    minimum size=3.5mm]
  \tikzstyle{v-p} = [circle, draw, fill=lightpurple,inner sep=0pt, 
    minimum size=3.5mm]
  \tikzstyle{cy2} = [draw]
  \tikzstyle{cy2-r} = [draw]
  \tikzstyle{cy2-s} = [draw,Blue]
  \tikzstyle{cy2-rs} = [draw,Green]
  \tikzstyle{cy2-r2s} = [draw,Purple]
  \tikzstyle{cy2-r3s} = [draw,orange] 
    \tikzstyle{cy2} = [draw]
  \tikzstyle{v-gr} = [circle, draw, fill=midgrey,inner sep=0pt,
    minimum size=5mm]
  \tikzstyle{v-big} = [circle, draw, fill=lightgrey,inner sep=0pt,
    minimum size=5mm]
      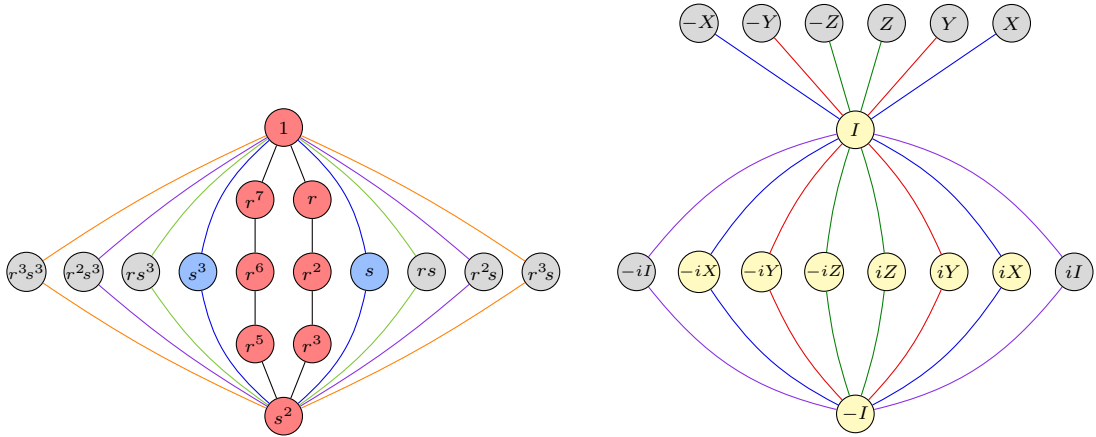
\begin{figure}[!ht]
  \begin{tikzpicture}[scale=1,yscale=.8,xscale=.95]
     \tikzstyle{v-small} = [circle, draw, fill=lightgrey,inner sep=0pt,
    minimum size=5mm]
    \begin{scope}[shift={(0,0)},scale=.8]
     \tikzstyle{every node}=[font=\scriptsize]
     \tikzstyle{v-r} = [circle, draw, fill=lred,inner sep=0pt,
    minimum size=5mm]
   \tikzstyle{v-b} = [circle, draw, fill=lblue,inner sep=0pt,
    minimum size=5mm]
    \tikzstyle{v-p} = [circle, draw, fill=lightpurple,inner sep=0pt, 
    minimum size=5mm]
      \node (1) at (0,3) [v-r] {$1$};
      \node (s2) at (0,-3) [v-r] {$s^2$};
      \node (s3) at (-1.5,0) [v-b] {$s^3$};
      \node (s) at (1.5,0) [v-b] {$s$};
      \node (rs3) at (-2.5,0) [v-small] {$rs^3$};
      \node (rs) at (2.5,0) [v-small] {$rs$};
      \node (r2s3) at (-3.5,0) [v-small] {\tiny $r^2\!s^3$};
      \node (r2s) at (3.5,0) [v-small] {$r^2\!s$};
      \node (r3s3) at (-4.5,0) [v-small] {\tiny $r^3\!s^3$};
      \node (r3s) at (4.5,0) [v-small] {$r^3\!s$};
      \draw [cy2-s,bend left=17] (1) to (s);
      \draw [cy2-s,bend right=17] (s2) to (s);
      \draw [cy2-rs,bend left=10] (1) to (rs);
      \draw [cy2-rs,bend right=10] (s2) to (rs);
      \draw [cy2-r2s,bend left=5] (1) to (r2s);
      \draw [cy2-r2s,bend right=5] (s2) to (r2s);
      \draw [cy2-r3s,bend left=5] (1) to (r3s); 
      \draw [cy2-r3s,bend right=5] (s2) to (r3s);
      \draw [cy2-s,bend right=17] (1) to (s3);
      \draw [cy2-s,bend left=17] (s2) to (s3);
      \draw [cy2-rs,bend right=10] (1) to (rs3);
      \draw [cy2-rs,bend left=10] (s2) to (rs3);
      \draw [cy2-r2s,bend right=5] (1) to (r2s3);
      \draw [cy2-r2s,bend left=5] (s2) to (r2s3);
      \draw [cy2-r3s,bend right=5] (1) to (r3s3); 
      \draw [cy2-r3s,bend left=5] (s2) to (r3s3);
      \node (r) at (.5,1.5) [v-r] {$r$};
      \node (r2) at (.5,0) [v-r] {$r^2$};
      \node (r3) at (.5,-1.5) [v-r] {$r^3$};
      \node (rs2) at (-.5,-1.5) [v-r] {$r^5$};
      \node (r2s2) at (-.5,0) [v-r] {$r^6$};
      \node (r3s2) at (-.5,1.5) [v-r] {$r^7$};

      \draw [cy2] (1) to (r);
      \draw [cy2] (r) to (r2);
      \draw [cy2] (r2) to (r3);
      \draw [cy2] (r3) to (s2);
      \draw [cy2] (s2) to (rs2);
      \draw [cy2] (rs2) to (r2s2);
      \draw [cy2] (r2s2) to (r3s2);
      \draw [cy2] (r3s2) to (1);
    \end{scope}
    \begin{scope}[shift={(8,0)},scale=1.75,yscale=1.35]
        \tikzstyle{v-yel} = [circle, draw, fill=lyellow,inner sep=0pt,
    minimum size=5mm]
      \tikzstyle{every node}=[font=\scriptsize]
      \node (-X) at (-1.25,1.75) [v-small] {$-X$};
      \node (-Y) at (-.75,1.75) [v-small] {$-Y$};
      \node (-Z) at (-.25,1.75) [v-small] {$-Z$};
      \node (X) at (1.25,1.75) [v-small] {$X$};
      \node (Y) at (.75,1.75) [v-small] {$Y$};
      \node (Z) at (.25,1.75) [v-small] {$Z$};      
      \node (I) at (0,1) [v-yel] {$I$};
      \node (-I) at (0,-1) [v-yel] {$-I$};
      \node (-iI) at (-1.75,0) [v-big] {\tiny$-iI$};
      \node (-iX) at (-1.25,0) [v-yel] {\tiny$-iX$};
      \node (-iY) at (-.75,0) [v-yel] {\tiny$-iY$};
      \node (-iZ) at (-.25,0) [v-yel] {\tiny$-iZ$};
      \node (iI) at (1.75,0) [v-big] {$iI$};
      \node (iX) at (1.25,0) [v-yel] {$iX$};
      \node (iY) at (.75,0) [v-yel] {$iY$};
      \node (iZ) at (.25,0) [v-yel] {$iZ$};
      \draw [cy2,Blue] (I) to (-X);
      \draw [cy2,Red] (I) to (-Y);
      \draw [cy2,darkgreen] (I) to (-Z);
      \draw [cy2,Blue] (I) to (X);
      \draw [cy2,Red] (I) to (Y);
      \draw [cy2,darkgreen] (I) to (Z);
      \draw [cy2,amethyst] (I) to [bend right=18] (-iI);
      \draw [cy2,Blue] (I) to [bend right=14] (-iX);
      \draw [cy2,Red] (I) to [bend right=10] (-iY);
      \draw [cy2,darkgreen] (I) to [bend right=6] (-iZ);
      \draw [cy2,amethyst] (I) to [bend left=18] (iI);
      \draw [cy2,Blue] (I) to [bend left=14] (iX);
      \draw [cy2,Red] (I) to [bend left=10] (iY);
      \draw [cy2,darkgreen] (I) to [bend left=6] (iZ);
      \draw [cy2,amethyst] (-I) to [bend right=18] (iI);
      \draw [cy2,Blue] (-I) to [bend right=14] (iX);
      \draw [cy2,Red] (-I) to [bend right=10] (iY);
      \draw [cy2,darkgreen] (-I) to [bend right=6] (iZ);
      \draw [cy2,amethyst] (-I) to [bend left=18] (-iI);
      \draw [cy2,Blue] (-I) to [bend left=14] (-iX);
      \draw [cy2,Red] (-I) to [bend left=10] (-iY);
      \draw [cy2,darkgreen] (-I) to [bend left=6] (-iZ);
    \end{scope}
    \end{tikzpicture}
  \caption{The cycle graphs for the generalized quaternion group $\Quat_{16}=\Dic_8$ (left) and the diquaternion group $\DQ_8$ (right). The yellow nodes in $\DQ_8$ correspond to the subgroup $\Quat_8$.}\label{fig:Q16-DQ8-cycle-graphs}
  \end{figure}
  
\begin{figure}[!ht]
  \begin{tikzpicture}[shorten >= -3pt, shorten <= -3pt, scale=1.6]
         \tikzstyle{every node}=[font=\small]
    \begin{scope}[shift={(0,0)}]
    \node (G) at (0,3) {$\DQ_8$};
    \node (C4xC2-1) at (-3,2) {$\Cyc_4\!\times\!\Cyc_2$}; 
    \node[faded] (D4-1) at (-2,2) {$\Dih_4$}; 
    \node (C4xC2-2) at (-1,2) {$\Cyc_4\!\times\!\Cyc_2$}; 
    \node[faded] (D4-2) at (0,2) {$\Dih_4$};  
    \node (C4xC2-3) at (1,2) {$\Cyc_4\!\times\!\Cyc_2$}; 
    \node[faded] (D4-3) at (2,2) {$\Dih_4\!=\!\bm{K}\!\!\!\!\!\!\!$};  
    \node[Red] (Q8) at (3,2) {$\Quat_8\!=\!\bm{H}\!\!\!\!\!\!\!\!$};  
    \node[faded] (V4-12) at (-3,1) {$V_4$};
    \node[faded] (V4-13) at (-2,1) {$V_4$};
    \node[faded] (V4-23) at (-1,1) {$V_4$};
    \node[Blue] (Z) at (0,1) {$\hspace{3mm}\Palert{\Cyc_4}\!=\!\bm{N}$};
    \node (C4-1) at (1,1) {$\Cyc_4$};
    \node (C4-2) at (2,1) {$\Cyc_4$};
    \node (C4-3) at (3,1) {$\Cyc_4$};
    \node[faded] (C2-11) at (-3.25,0) {$\Cyc_2$};
    \node[faded] (C2-12) at (-2.75,0) {$\Cyc_2$};
    \node[faded] (C2-21) at (-2.25,0) {$\Cyc_2$};
    \node[faded] (C2-22) at (-1.75,0) {$\Cyc_2$};
    \node[faded] (C2-31) at (-1.25,0) {$\Cyc_2$};
    \node[faded] (C2-32) at (-0.75,0) {$\Cyc_2$};    
    \node[amethyst] (-I) at (1,0) {$\Cyc_2\!=\!\bm{N\cap H\!=\!N\cap K}$\!\!\!\!\!\!\!\!\!\!\!\!\!\!\!\!\!\!\!\!\!\!\!\!\!\!\!\!\!\!\!\!\!\!\!\!\!\!\!\!\!\!\!\!\!\!\!\!\!\!};
    \node[amethyst] (I) at (0,-1) {$\Cyc_1$};
    \draw (C4xC2-1) to (G); 
    \draw[faded] (D4-1) to (G); 
    \draw (C4xC2-2) to (G);
    \draw[faded] (D4-2) to (G); 
    \draw (C4xC2-3) to (G); 
    \draw[faded] (D4-3) to (G);
    \draw (Q8) to (G);
    \draw (C4xC2-1) to (Z); 
    \draw (C4xC2-2) to (Z); 
    \draw (C4xC2-3) to (Z);
    \draw[faded] (V4-12) to (D4-1); 
    \draw[faded] (V4-12) to (D4-2); 
    \draw[faded] (V4-12) to (C4xC2-3);
    \draw[faded] (V4-13) to (D4-1); 
    \draw[faded] (V4-13) to (D4-3); 
    \draw[faded] (V4-13) to (C4xC2-2);
    \draw[faded] (V4-23) to (D4-2); 
    \draw[faded] (V4-23) to (D4-3); 
    \draw[faded] (V4-23) to (C4xC2-1);
    \draw[faded] (C4-1) to (C4xC2-1); \draw[faded] (C4-1) to (D4-1);
    \draw[faded] (C4-2) to (C4xC2-2); \draw[faded] (C4-2) to (D4-2);
    \draw[faded] (C4-3) to (C4xC2-3); \draw[faded] (C4-3) to (D4-3);
    \draw (Q8) to (C4-1); 
    \draw (Q8) to (C4-2); 
    \draw (Q8) to (C4-3);
    \draw[faded] (V4-12) to (C2-11); \draw[faded] (V4-12) to (C2-12); 
    \draw[faded] (V4-13) to (C2-21); \draw[faded] (V4-13) to (C2-22); 
    \draw[faded] (V4-23) to (C2-31); \draw[faded] (V4-23) to (C2-32); 
    \draw[faded] (-I) to (V4-12); 
    \draw[faded] (-I) to (V4-13); 
    \draw[faded] (-I) to (V4-23);
    \draw (-I) to (C4-1); 
    \draw (-I) to (C4-2); 
    \draw (-I) to (C4-3);
    \draw[amethyst,very thick] (-I) to (Z);
    \draw[faded] (C2-11) to (I); \draw[faded] (C2-12) to (I); 
    \draw[faded] (C2-21) to (I); \draw[faded] (C2-22) to (I); 
    \draw[faded] (C2-31) to (I); \draw[faded] (C2-32) to (I);
    \draw[amethyst,very thick] (-I) to (I);
  \end{scope}
  \end{tikzpicture}
\caption{From the subgroup lattice of the diquaternion group, we can see that it is the central product $\Quat_8\circ \Cyc_4$ of $H=\<XY,XZ\>$ with $N=\<-iI\>$, and a central product $\Dih_4\circ \Cyc_4$, where $K=\<X,Y\>\cong \Dih_4$.}\label{fig:DQ8-central-product}
\end{figure}
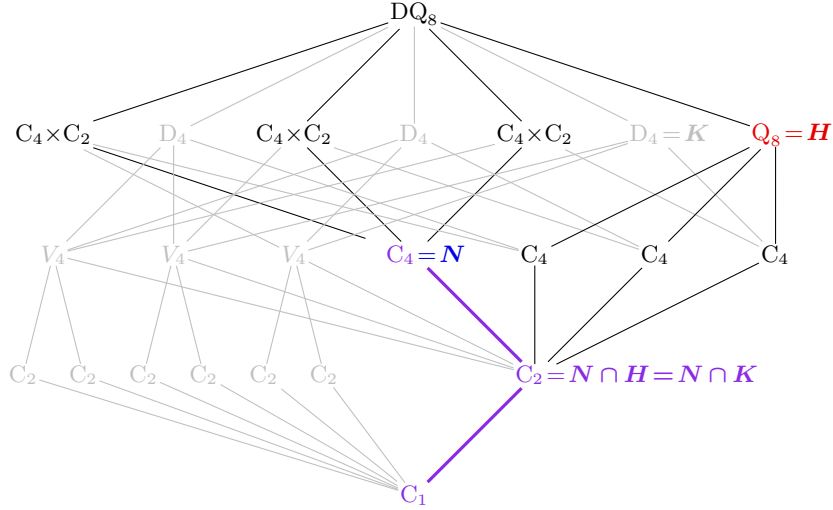

  The subgroup lattice of the diquaternion group, shown in Figure~\ref{fig:DQ8-central-product}, also reveals information about its structure. If we ``chop off'' the lattice at the normal subgroup $\<-I\>\cong \Cyc_2$, it is apparent that the order-$8$ quotient $\DQ_8/\<-I\>$ has seven order-$2$ subgroups, so it must be isomorphic to the abelian group $\Cyc_2^3$. Moreover, $\DQ_8$ is the product of the subgroup $H=\<XY,XZ\>\cong \Quat_8$ and the unicorn $N=\<-iI\>\cong \Cyc_4$ --- the only order-$4$ subgroup contained in three abelian subgroups. Since $N$ is central (it actually is the center), $\DQ_8$ is a \emph{central product}, $\Quat_8\circ \Cyc_4$. Similarly, it is also a central product of $K=\<X,Y\>\cong \Dih_4$ with $N$.

Our ``dihedralizing the quaternions'' construction --- adding a reflection from $\Dih_n$ to $\Quat_8$, works for \emph{any} generalized quaternion group, $\Quat_n$ (recall, $n=2^m$). Appropriately, we will call such a group a \emph{generalized diquaternion group}. Alternatively, we can think of this as the result of starting with the diquaternion group $\DQ_8=\<i,j,f\>$, and replacing $i=\zeta_4=e^{2\pi i/4}$ with a larger root of unity, $\zeta_n=e^{2\pi i/n}$, like we did to construct the dicyclic groups from $\Quat_8$. That is, define
   \[
  \DQ_n:=\big\<\zeta_n,j,\zeta_nj,f\big\>\cong\left\<
  \underbrace{\begin{bmatrix} \zeta_n & 0 \\ 0 & \bar{\zeta}_n\end{bmatrix}}_{R=R_n},\;
  \underbrace{\begin{bmatrix} 0 & -1 \\ 1 & 0\end{bmatrix}}_{S},\;
  \underbrace{\begin{bmatrix} 0 & -\zeta_n \\ \bar\zeta_n & 0\end{bmatrix}}_{T=T_n},\;
  \Balert{\underbrace{\begin{bmatrix}0 & 1 \\ 1 & 0\end{bmatrix}}_{\Balert{F}}}\right\>.
  \]
 This group can also be constructed by appropriately modifying the Pauli matrices. Since $X$ and $Z$ have only $\pm 1$ entries, they are unchanged. The entries in $Y_4:=Y$, which are  $i$ and $-i$, are replaced with $\zeta_n$ and $\bar\zeta_n$ to get $Y_n$. A Cayley graph of $\DQ_{16}$ is shown in Figure~\ref{fig:DQ16-Cayley}. Once again, the highlighted nodes denote the generalized quaternion subgroup, $\Quat_{16}$. We can also see an index-$2$ dihedral subgroup $\<X,Y\>\cong \Dih_8$ from the nodes along the outer ring, and an index-$2$ abelian subgroup $\<XY,Z\>\cong \Cyc_8\times \Cyc_2$ from the nodes at angles of $2\pi ik/4$ radians, for some $k\in\Z$.

\tikzstyle{v-small} = [circle, draw, fill=lightgrey,inner sep=0pt, 
  minimum size=4mm]
\tikzstyle{v-yel} = [circle, draw, fill=lyellow,inner sep=0pt, 
  minimum size=4mm]
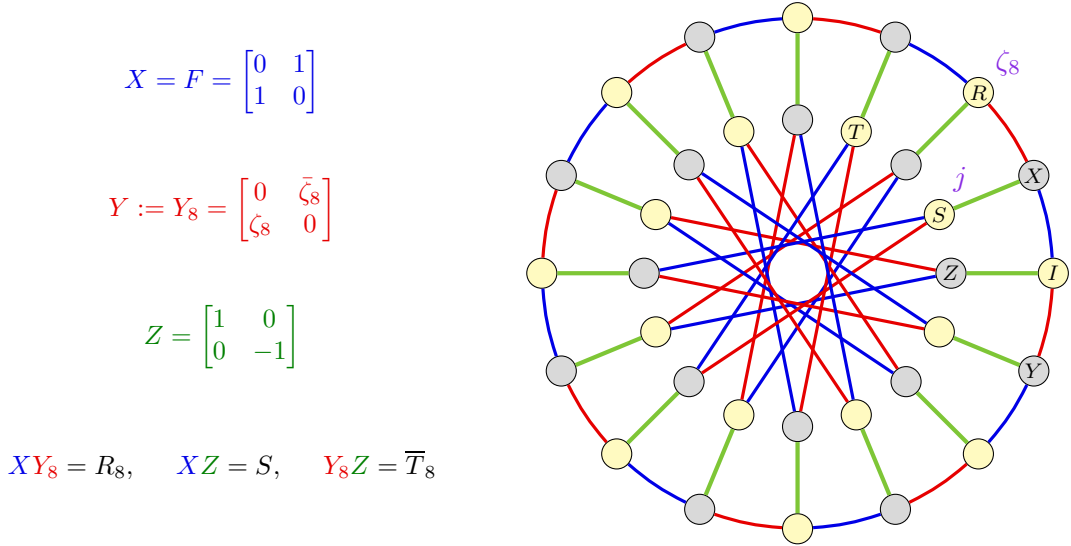
\begin{figure}[!ht]
\begin{tikzpicture}[scale=1.7]
  \tikzstyle{every node}=[font=\footnotesize]
  \tikzstyle{rr-bend} = [draw, very thick, Red,bend left=7]
  \tikzstyle{bb-bend} = [draw, very thick, Blue,bend left=7]
  \tikzstyle{rr} = [draw, very thick, Red]
  \tikzstyle{bb} = [draw, very thick, Blue]
  \tikzstyle{gg} = [draw, ultra thick, Green]
  \begin{scope}[shift={(0,0)}]
    \node at (-4.5,1.5)
          {\normalsize\Balert{$X=F=\begin{bmatrix}0&1\\1&0\end{bmatrix}$}};
    \node at (-4.5,.5)
          {\normalsize\Alert{$Y:=Y_8=\begin{bmatrix}0&\bar{\zeta}_8\\ \zeta_8 &0\end{bmatrix}$}};
    \node at (-4.5,-.5)
          {\normalsize\Galert{$Z=\begin{bmatrix}1&0\\0&-1\end{bmatrix}$}};
          \node at (-4.5,-1.5) {\normalsize $\Balert{X}\Alert{Y_8}=R_8$,\quad\; $\Balert{X}\Galert{Z}=S$,\quad\; $\Alert{Y_8}\Galert{Z}=\overline{T}_8$};
      \node at (45:2.33) {\large\Palert{$\zeta_8$}};
      \node at (30:1.47) {\large\Palert{$j$}};
      \node (s) at (90:1.2) [v-small] {};
      \node (rs) at (67.5:1.2) [v-yel] {$T$};
      \node (r2s) at (45:1.2) [v-small] {};
      \node (r3s) at (22.5:1.2) [v-yel] {$S$};
      \node (r4s) at (0:1.2) [v-small] {$Z$};
      \node (r5s) at (-22.5:1.2) [v-yel] {};
      \node (r6s) at (-45:1.2) [v-small] {};
      \node (r7s) at (-67.5:1.2) [v-yel] {};
      \node (r8s) at (-90:1.2) [v-small] {};
      \node (r9s) at (-112.5:1.2) [v-yel] {};
      \node (r10s) at (-135:1.2) [v-small] {};
      \node (r11s) at (-157.5:1.2) [v-yel] {};
      \node (r12s) at (180:1.2) [v-small] {};
      \node (r13s) at (157.5:1.2) [v-yel] {};
      \node (r14s) at (135:1.2) [v-small] {};
      \node (r15s) at (112.5:1.2) [v-yel] {};
      \node (1) at (90:2) [v-yel] {};
      \node (r) at (67.5:2) [v-small] {};
      \node (r2) at (45:2) [v-yel] {$R$};
      \node (r3) at (22.5:2) [v-small] {$X$};
      \node (r4) at (0:2) [v-yel] {$I$};
      \node (r5) at (-22.5:2) [v-small] {$Y$};
      \node (r6) at (-45:2) [v-yel] {};
      \node (r7) at (-67.5:2) [v-small] {};
      \node (r8) at (-90:2) [v-yel] {};
      \node (r9) at (-112.5:2) [v-small] {};
      \node (r10) at (-135:2) [v-yel] {};
      \node (r11) at (-157.5:2) [v-small] {};
      \node (r12) at (180:2) [v-yel] {};
      \node (r13) at (157.5:2) [v-small] {};
      \node (r14) at (135:2) [v-yel] {};
      \node (r15) at (112.5:2) [v-small] {};
      \draw [rr-bend] (1) to (r); \draw [bb-bend] (r) to (r2);
      \draw [rr-bend] (r2) to (r3); \draw [bb-bend] (r3) to (r4);
      \draw [rr-bend] (r4) to (r5); \draw [bb-bend] (r5) to (r6);
      \draw [rr-bend] (r6) to (r7); \draw [bb-bend] (r7) to (r8);
      \draw [rr-bend] (r8) to (r9); \draw [bb-bend] (r9) to (r10);
      \draw [rr-bend] (r10) to (r11); \draw [bb-bend] (r11) to (r12);
      \draw [rr-bend] (r12) to (r13); \draw [bb-bend] (r13) to (r14);
      \draw [rr-bend] (r14) to (r15); \draw [bb-bend] (r15) to (1);
      \draw [rr] (s) to (r9s); \draw [bb] (r9s) to (r2s);
      \draw [rr] (r2s) to (r11s); \draw [bb] (r11s) to (r4s);
      \draw [rr] (r4s) to (r13s); \draw [bb] (r13s) to (r6s);
      \draw [rr] (r6s) to (r15s); \draw [bb] (r15s) to (r8s);
      \draw [rr] (r8s) to (rs); \draw [bb] (rs) to (r10s);
      \draw [rr] (r10s) to (r3s); \draw [bb] (r3s) to (r12s);
      \draw [rr] (r12s) to (r5s); \draw [bb] (r5s) to (r14s);
      \draw [rr] (r14s) to (r7s); \draw [bb] (r7s) to (s); 
      \draw [gg] (1) to (s); \draw [gg] (r) to (rs);
      \draw [gg] (r2) to (r2s); \draw [gg] (r3) to (r3s);
      \draw [gg] (r4) to (r4s); \draw [gg] (r5) to (r5s);
      \draw [gg] (r6) to (r6s); \draw [gg] (r7) to (r7s);
      \draw [gg] (r8) to (r8s); \draw [gg] (r9) to (r9s);
      \draw [gg] (r10) to (r10s); \draw [gg] (r11) to (r11s);
      \draw [gg] (r12) to (r12s); \draw [gg] (r13) to (r13s);
      \draw [gg] (r14) to (r14s); \draw [gg] (r15) to (r15s);
  \end{scope}
  \end{tikzpicture}
  \caption{A Cayley graph of the generalized diquaternion group $\DQ_{16}$, generated by the ``generalized Pauli matrices.'' The highlighted nodes show the generalized quaternion group $\Quat_{16}$ as a subgroup.}\label{fig:DQ16-Cayley}
  \end{figure}
  
  The generalized diquaternion group $\DQ_{16}$ is a central product of both the subgroups $\Dih_8$ and $\Quat_{16}$ with the central cyclic subgroup $\<-iI\>\cong \Cyc_4$. It is worth noting that the Pauli group on $1$ qubit generalizes to the Pauli group on $n$ qubits by tensoring the matrices. However, not all generalized diquaternion groups are Pauli groups. There is a diquaternion group of order $n=2^m$ for all $m\geq 4$, whereas the Pauli group on $n$ qubits has order $4^{n+1}$. To the best of our knowledge, the family of diquaternion groups does not exist in the literature with an alternate name.

\section{Semidihedral and semiabelian groups} 

  The dihedral group $\Dih_n$ contains a normal index-$2$ cyclic subgroup $\<r\>$, and an involution $f\not\in\<r\>$. Therefore, it is the semidirect product of $\Cyc_n=\<r\>$ and $\Cyc_2=\<f\>$. It is natural to ask how many other groups there are with an index-$2$ cyclic subgroup. Of course, $\Cyc_n\times \Cyc_2$ is another example. Any group that is a semidirect product of $\Cyc_n=\<r\>$ with $\Cyc_2$ must have an involution $s\not\in\<r\>$. In terms of the Cayley graph, we are asking how many ways there are to complete a ``partial Cayley graph'' like the one shown in the middle of Figure~\ref{fig:partial-Cayley}. For sake of comparison, the dicyclic groups have an index-$2$ subgroup, but all other elements have order $4$. Thus, one can think about their construction as a way to ``re-wire'' the \emph{blue} $f$-arrows in a Cayley graph for $\Dih_n=\<r,f\>$, as in Figure~\ref{fig:partial-Cayley} (right). In contrast, our question about semidirect products is asking: \emph{how can we rewire the inner red $r$-arrows in the dihedral group?}

  \tikzstyle{v-small} = [circle, draw, fill=lightgrey,inner sep=0pt, 
    minimum size=4.5mm]
  \tikzstyle{R} = [draw, very thick, Red,-stealth,bend right=10]
  \tikzstyle{R2} = [draw, very thick, Red,-stealth,bend right=15]
  \tikzstyle{R3} = [draw, very thick, Red,-stealth,bend left=12]
  \tikzstyle{B} = [draw, very thick, Blue,-stealth,bend right=25]
  \begin{figure}[!ht]
  \begin{tikzpicture}[scale=1,auto]
    \tikzstyle{every node}=[font=\tiny]
    \begin{scope}[shift={(0,0)}]
      \node at (0,0) {\normalsize $\Dih_8$};
      \node (s) at (0:1.1) [v-small] {$s$};
      \node (rs) at (45:1.15) [v-small] {$rs$};
      \node (r2s) at (90:1.15) [v-small] {$r^2\!s$};
      \node (r3s) at (135:1.15) [v-small] {$r^3\!s$};
      \node (r4s) at (180:1.15) [v-small] {$r^4\!s$};
      \node (r5s) at (225:1.15) [v-small] {$r^5\!s$};
      \node (r6s) at (270:1.15) [v-small] {$r^6\!s$};
      \node (r7s) at (315:1.15) [v-small] {$r^7\!s$};
      \node (1) at (0:2) [v-small] {$1$};
      \node (r) at (45:2) [v-small] {$r$};
      \node (r2) at (90:2) [v-small] {$r^2$};
      \node (r3) at (135:2) [v-small] {$r^3$};
      \node (r4) at (180:2) [v-small] {$r^4$};
      \node (r5) at (225:2) [v-small] {$r^5$};
      \node (r6) at (270:2) [v-small] {$r^6$};
      \node (r7) at (315:2) [v-small] {$r^7$};
      \draw [R2] (1) to (r);
      \draw [R2] (r) to (r2);
      \draw [R2] (r2) to (r3);
      \draw [R2] (r3) to (r4);
      \draw [R2] (r4) to (r5);
      \draw [R2] (r5) to (r6);
      \draw [R2] (r6) to (r7);
      \draw [R2] (r7) to (1);
      \draw [R3] (s) to (r7s);
      \draw [R3] (r7s) to (r6s);
      \draw [R3] (r6s) to (r5s);
      \draw [R3] (r5s) to (r4s);
      \draw [R3] (r4s) to (r3s);
      \draw [R3] (r3s) to (r2s);
      \draw [R3] (r2s) to (rs);
      \draw [R3] (rs) to (s);
      \draw [bb] (1) to (s); \draw [bb] (r) to (rs);
      \draw [bb] (r2) to (r2s); \draw [bb] (r3) to (r3s);
      \draw [bb] (r4) to (r4s); \draw [bb] (r5) to (r5s);
      \draw [bb] (r6) to (r6s); \draw [bb] (r7) to (r7s);
    \end{scope}
    \begin{scope}[shift={(10,0)}]
      \node (s) at (0:1.1) [v-small] {$s$};
      \node (rs) at (45:1.15) [v-small] {$rs$};
      \node (r2s) at (90:1.15) [v-small] {$r^2\!s$};
      \node (r3s) at (135:1.15) [v-small] {$r^3\!s$};
      \node (r4s) at (180:1.15) [v-small] {$r^4\!s$};
      \node (r5s) at (225:1.15) [v-small] {$r^5\!s$};
      \node (r6s) at (270:1.15) [v-small] {$r^6\!s$};
      \node (r7s) at (315:1.15) [v-small] {$r^7\!s$};
      \node (1) at (0:2) [v-small] {$1$};
      \node (r) at (45:2) [v-small] {$r$};
      \node (r2) at (90:2) [v-small] {$r^2$};
      \node (r3) at (135:2) [v-small] {$r^3$};
      \node (r4) at (180:2) [v-small] {$r^4$};
      \node (r5) at (225:2) [v-small] {$r^5$};
      \node (r6) at (270:2) [v-small] {$r^6$};
      \node (r7) at (315:2) [v-small] {$r^7$};
      \draw [R2] (1) to (r);
      \draw [R2] (r) to (r2);
      \draw [R2] (r2) to (r3);
      \draw [R2] (r3) to (r4);
      \draw [R2] (r4) to (r5);
      \draw [R2] (r5) to (r6);
      \draw [R2] (r6) to (r7);
      \draw [R2] (r7) to (1);
      \draw [R3] (s) to (r7s);
      \draw [R3] (r7s) to (r6s);
      \draw [R3] (r6s) to (r5s);
      \draw [R3] (r5s) to (r4s);
      \draw [R3] (r4s) to (r3s);
      \draw [R3] (r3s) to (r2s);
      \draw [R3] (r2s) to (rs);
      \draw [R3] (rs) to (s);
    \end{scope}
    \begin{scope}[shift={(5,0)}]
      \node (s) at (0:1.1) [v-small] {$s$};
      \node (rs) at (45:1.15) [v-small] {$rs$};
      \node (r2s) at (90:1.15) [v-small] {$r^2\!s$};
      \node (r3s) at (135:1.15) [v-small] {$r^3\!s$};
      \node (r4s) at (180:1.15) [v-small] {$r^4\!s$};
      \node (r5s) at (225:1.15) [v-small] {$r^5\!s$};
      \node (r6s) at (270:1.15) [v-small] {$r^6\!s$};
      \node (r7s) at (315:1.15) [v-small] {$r^7\!s$};
      \node (1) at (0:2) [v-small] {$1$};
      \node (r) at (45:2) [v-small] {$r$};
      \node (r2) at (90:2) [v-small] {$r^2$};
      \node (r3) at (135:2) [v-small] {$r^3$};
      \node (r4) at (180:2) [v-small] {$r^4$};
      \node (r5) at (225:2) [v-small] {$r^5$};
      \node (r6) at (270:2) [v-small] {$r^6$};
      \node (r7) at (315:2) [v-small] {$r^7$};
      \draw [R2] (1) to (r);
      \draw [R2] (r) to (r2);
      \draw [R2] (r2) to (r3);
      \draw [R2] (r3) to (r4);
      \draw [R2] (r4) to (r5);
      \draw [R2] (r5) to (r6);
      \draw [R2] (r6) to (r7);
      \draw [R2] (r7) to (1);
      \draw [bb] (1) to (s); \draw [bb] (r) to (rs);
      \draw [bb] (r2) to (r2s); \draw [bb] (r3) to (r3s);
      \draw [bb] (r4) to (r4s); \draw [bb] (r5) to (r5s);
      \draw [bb] (r6) to (r6s); \draw [bb] (r7) to (r7s);
    \end{scope}
  \end{tikzpicture}
  \caption{Any semidirect product of $\Cyc_8=\<r\>$ with $\Cyc_2=\<s\>$ must have a ``partial Cayley graph'' as shown in the middle. The group $\Quat_{16}=\Dic_8$ is not a semidirect product, but it can be built from the partial Cayley graph shown on the right.}\label{fig:partial-Cayley}
  \end{figure}
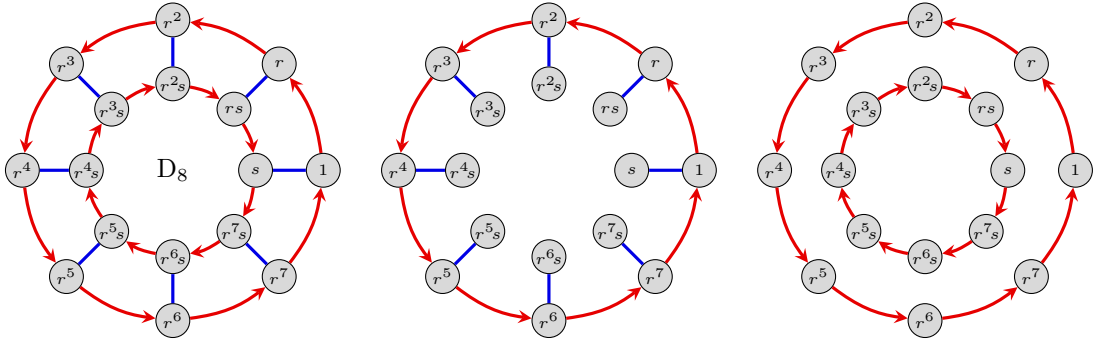

  A semidirect product $\Cyc_n\rtimes_\theta\Cyc_2$ of cyclic groups is defined by an automorphism $\theta\colon \Cyc_2\to\Aut(\Cyc_n)$. In this section, we will learn about two relatively obscure families of groups that arise when $n=2^m$. The automorphism group of $\Cyc_n$ is
  \[
  \Aut(\Cyc_n)=\big\{\sigma_k\mid\;1\leq k<n,\; \gcd(n,k)=1\big\},\qquad\text{where $\sigma_k(1)=k$},
  \]
  and is isomorphic to the group $U_n$ of units of $\Z_n$. Any automorphism $\Cyc_2\to\Aut(\Cyc_{2^n})$ must send $1$ to a generator of $U_{2^n}$. It is a basic number theory fact that $x^2\equiv 1$ has four solutions modulo $2^n$, which are $\pm 1$ and $2^{n-1}\!\pm\!1$. Thus, there are exactly four semidirect products of $\Cyc_{2^n}$ with $\Cyc_2$. All four are generated by $r$ and $s$, subject to the relations $r^{2^n}=1$ and $s^2=1$. The only difference is the last relation, which is of the form $srs=r^k$. The four possibilities are $k=1$ (the abelian group $\Cyc_{2^n}\times \Cyc_2$), $k=-1=2^n-1$ (the dihedral group $\Dih_{2^n}$), $k=2^{n-1}-1$ (the \emph{semidihedral group} $\SD_{2^n}$), and $k=2^{n-1}+1$ (the \emph{semiabelian group} $\SA_{2^n}$). Partial Cayley graphs of these for $2^n=16$ are shown in Figure~\ref{fig:semidirect-products}, which are meant to highlight that last relation.

  \tikzstyle{v-small} = [circle, draw, fill=lightgrey,inner sep=0pt, 
    minimum size=4.5mm]
    \tikzstyle{vv} = [circle, draw, lightgrey, fill=lightgrey,inner sep=0pt, 
    minimum size=1mm]
     \tikzstyle{VV} = [circle, draw, fill=lightgrey,inner sep=0pt, 
    minimum size=1mm]
 
  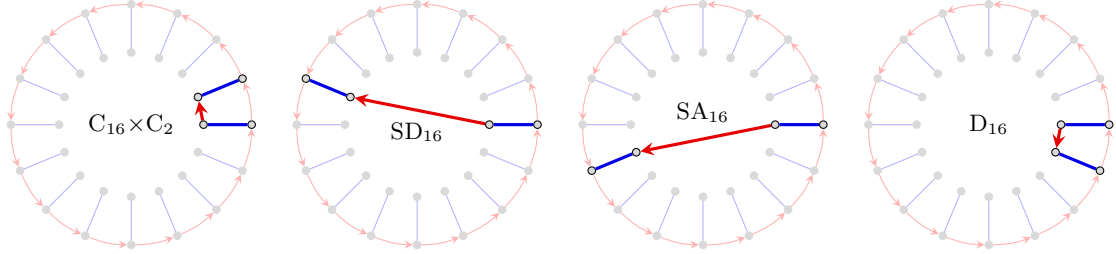
\begin{figure}[!ht]
   \tikzstyle{R} = [draw, rfaded,-stealth,bend right=10]
   \tikzstyle{BB} = [draw, bfaded]
  \begin{tikzpicture}[scale=.8,auto]
    \tikzstyle{every node}=[font=\small]
    \begin{scope}[shift={(0,0)}]
      \node (s) at (0:1.2) [VV] {};
      \node (rs) at (22.5:1.2) [VV] {};
      \node (r2s) at (45:1.2) [vv] {};
      \node (r3s) at (67.5:1.2) [vv] {};
      \node (r4s) at (90:1.2) [vv] {};
      \node (r5s) at (112.5:1.2) [vv] {};
      \node (r6s) at (135:1.2) [vv] {};
      \node (r7s) at (157.5:1.2) [vv] {};
      \node (r8s) at (180:1.2) [vv] {};
      \node (r9s) at (-157.5:1.2) [vv] {};
      \node (r10s) at (-135:1.2) [vv] {};
      \node (r11s) at (-112.5:1.2) [vv] {};
      \node (r12s) at (-90:1.2) [vv] {};
      \node (r13s) at (-67.5:1.2) [vv] {};
      \node (r14s) at (-45:1.2) [vv] {};
      \node (r15s) at (-22.5:1.2) [vv] {};
      \node (1) at (0:2) [VV] {};
      \node (r) at (22.5:2) [VV] {};
      \node (r2) at (45:2) [vv] {};
      \node (r3) at (67.5:2) [vv] {};
      \node (r4) at (90:2) [vv] {};
      \node (r5) at (112.5:2) [vv] {};
      \node (r6) at (135:2) [vv] {};
      \node (r7) at (157.5:2) [vv] {};
      \node (r8) at (180:2) [vv] {};
      \node (r9) at (-157.5:2) [vv] {};
      \node (r10) at (-135:2) [vv] {};
      \node (r11) at (-112.5:2) [vv] {};
      \node (r12) at (-90:2) [vv] {};
      \node (r13) at (-67.5:2) [vv] {};
      \node (r14) at (-45:2) [vv] {};
      \node (r15) at (-22.5:2) [vv] {};
      \draw [R] (1) to (r); \draw [R] (r) to (r2); \draw [R] (r2) to (r3);
      \draw [R] (r3) to (r4); \draw [R] (r4) to (r5); \draw [R] (r5) to (r6);
      \draw [R] (r6) to (r7); \draw [R] (r7) to (r8); \draw [R] (r8) to (r9);
      \draw [R] (r9) to (r10); \draw [R] (r10) to (r11);
      \draw [R] (r11) to (r12); \draw [R] (r12) to (r13);
      \draw [R] (r13) to (r14); \draw [R] (r14) to (r15);
      \draw [R] (r15) to (1);
      \draw [r] (s) to (rs);  
      \draw [bb] (1) to (s); \draw [bb] (r) to (rs);
      \draw [BB] (r2) to (r2s); \draw [BB] (r3) to (r3s);
      \draw [BB] (r4) to (r4s); \draw [BB] (r5) to (r5s);
      \draw [BB] (r6) to (r6s); \draw [BB] (r7) to (r7s);
      \draw [BB] (r8) to (r8s); \draw [BB] (r9) to (r9s);
      \draw [BB] (r10) to (r10s); \draw [BB] (r11) to (r11s);
      \draw [BB] (r12) to (r12s); \draw [BB] (r13) to (r13s);
      \draw [BB] (r14) to (r14s); \draw [BB] (r15) to (r15s);
      \node at (0,0) {$\Cyc_{16}\!\times\!\Cyc_2$};
    \end{scope}
    \begin{scope}[shift={(4.75,0)}]
      \node (s) at (0:1.2) [VV] {};
      \node (rs) at (22.5:1.2) [vv] {};
      \node (r2s) at (45:1.2) [vv] {};
      \node (r3s) at (67.5:1.2) [vv] {};
      \node (r4s) at (90:1.2) [vv] {};
      \node (r5s) at (112.5:1.2) [vv] {};
      \node (r6s) at (135:1.2) [vv] {};
      \node (r7s) at (157.5:1.2) [VV] {};
      \node (r8s) at (180:1.2) [vv] {};
      \node (r9s) at (-157.5:1.2) [vv] {};
      \node (r10s) at (-135:1.2) [vv] {};
      \node (r11s) at (-112.5:1.2) [vv] {};
      \node (r12s) at (-90:1.2) [vv] {};
      \node (r13s) at (-67.5:1.2) [vv] {};
      \node (r14s) at (-45:1.2) [vv] {};
      \node (r15s) at (-22.5:1.2) [vv] {};
      \node (1) at (0:2) [VV] {};
      \node (r) at (22.5:2) [vv] {};
      \node (r2) at (45:2) [vv] {};
      \node (r3) at (67.5:2) [vv] {};
      \node (r4) at (90:2) [vv] {};
      \node (r5) at (112.5:2) [vv] {};
      \node (r6) at (135:2) [vv] {};
      \node (r7) at (157.5:2) [VV] {};
      \node (r8) at (180:2) [vv] {};
      \node (r9) at (-157.5:2) [vv] {};
      \node (r10) at (-135:2) [vv] {};
      \node (r11) at (-112.5:2) [vv] {};
      \node (r12) at (-90:2) [vv] {};
      \node (r13) at (-67.5:2) [vv] {};
      \node (r14) at (-45:2) [vv] {};
      \node (r15) at (-22.5:2) [vv] {};
      \draw [R] (1) to (r); \draw [R] (r) to (r2); \draw [R] (r2) to (r3);
      \draw [R] (r3) to (r4); \draw [R] (r4) to (r5); \draw [R] (r5) to (r6);
      \draw [R] (r6) to (r7); \draw [R] (r7) to (r8); \draw [R] (r8) to (r9);
      \draw [R] (r9) to (r10); \draw [R] (r10) to (r11);
      \draw [R] (r11) to (r12); \draw [R] (r12) to (r13);
      \draw [R] (r13) to (r14); \draw [R] (r14) to (r15);
      \draw [R] (r15) to (1);
      \draw [r] (s) to (r7s);  
      \draw [bb] (1) to (s); \draw [BB] (r) to (rs);
      \draw [BB] (r2) to (r2s); \draw [BB] (r3) to (r3s);
      \draw [BB] (r4) to (r4s); \draw [BB] (r5) to (r5s);
      \draw [BB] (r6) to (r6s); \draw [bb] (r7) to (r7s);
      \draw [BB] (r8) to (r8s); \draw [BB] (r9) to (r9s);
      \draw [BB] (r10) to (r10s); \draw [BB] (r11) to (r11s);
      \draw [BB] (r12) to (r12s); \draw [BB] (r13) to (r13s);
      \draw [BB] (r14) to (r14s); \draw [BB] (r15) to (r15s);
      \node at (270:.12) {$\SD_{16}$};
    \end{scope}
    \begin{scope}[shift={(9.5,0)}]
      \node (s) at (0:1.2) [VV] {};
      \node (rs) at (22.5:1.2) [vv] {};
      \node (r2s) at (45:1.2) [vv] {};
      \node (r3s) at (67.5:1.2) [vv] {};
      \node (r4s) at (90:1.2) [vv] {};
      \node (r5s) at (112.5:1.2) [vv] {};
      \node (r6s) at (135:1.2) [vv] {};
      \node (r7s) at (157.5:1.2) [vv] {};
      \node (r8s) at (180:1.2) [vv] {};
      \node (r9s) at (-157.5:1.2) [VV] {};
      \node (r10s) at (-135:1.2) [vv] {};
      \node (r11s) at (-112.5:1.2) [vv] {};
      \node (r12s) at (-90:1.2) [vv] {};
      \node (r13s) at (-67.5:1.2) [vv] {};
      \node (r14s) at (-45:1.2) [vv] {};
      \node (r15s) at (-22.5:1.2) [vv] {};
      \node (1) at (0:2) [VV] {};
      \node (r) at (22.5:2) [vv] {};
      \node (r2) at (45:2) [vv] {};
      \node (r3) at (67.5:2) [vv] {};
      \node (r4) at (90:2) [vv] {};
      \node (r5) at (112.5:2) [vv] {};
      \node (r6) at (135:2) [vv] {};
      \node (r7) at (157.5:2) [vv] {};
      \node (r8) at (180:2) [vv] {};
      \node (r9) at (-157.5:2) [VV] {};
      \node (r10) at (-135:2) [vv] {};
      \node (r11) at (-112.5:2) [vv] {};
      \node (r12) at (-90:2) [vv] {};
      \node (r13) at (-67.5:2) [vv] {};
      \node (r14) at (-45:2) [vv] {};
      \node (r15) at (-22.5:2) [vv] {};
      \draw [R] (1) to (r); \draw [R] (r) to (r2); \draw [R] (r2) to (r3);
      \draw [R] (r3) to (r4); \draw [R] (r4) to (r5); \draw [R] (r5) to (r6);
      \draw [R] (r6) to (r7); \draw [R] (r7) to (r8); \draw [R] (r8) to (r9);
      \draw [R] (r9) to (r10); \draw [R] (r10) to (r11);
      \draw [R] (r11) to (r12); \draw [R] (r12) to (r13);
      \draw [R] (r13) to (r14); \draw [R] (r14) to (r15);
      \draw [R] (r15) to (1);
      \draw [r] (s) to (r9s);  
      \draw [bb] (1) to (s); \draw [BB] (r) to (rs);
      \draw [BB] (r2) to (r2s); \draw [BB] (r3) to (r3s);
      \draw [BB] (r4) to (r4s); \draw [BB] (r5) to (r5s);
      \draw [BB] (r6) to (r6s); \draw [BB] (r7) to (r7s);
      \draw [BB] (r8) to (r8s); \draw [bb] (r9) to (r9s);
      \draw [BB] (r10) to (r10s); \draw [BB] (r11) to (r11s);
      \draw [BB] (r12) to (r12s); \draw [BB] (r13) to (r13s);
      \draw [BB] (r14) to (r14s); \draw [BB] (r15) to (r15s);
      \node at (90:.2) {$\SA_{16}$};
    \end{scope}
    \begin{scope}[shift={(14.25,0)}]
      \node (s) at (0:1.2) [VV] {};
      \node (rs) at (22.5:1.2) [vv] {};
      \node (r2s) at (45:1.2) [vv] {};
      \node (r3s) at (67.5:1.2) [vv] {};
      \node (r4s) at (90:1.2) [vv] {};
      \node (r5s) at (112.5:1.2) [vv] {};
      \node (r6s) at (135:1.2) [vv] {};
      \node (r7s) at (157.5:1.2) [vv] {};
      \node (r8s) at (180:1.2) [vv] {};
      \node (r9s) at (-157.5:1.2) [vv] {};
      \node (r10s) at (-135:1.2) [vv] {};
      \node (r11s) at (-112.5:1.2) [vv] {};
      \node (r12s) at (-90:1.2) [vv] {};
      \node (r13s) at (-67.5:1.2) [vv] {};
      \node (r14s) at (-45:1.2) [vv] {};
      \node (r15s) at (-22.5:1.2) [VV] {};
      \node (1) at (0:2) [VV] {};
      \node (r) at (22.5:2) [vv] {};
      \node (r2) at (45:2) [vv] {};
      \node (r3) at (67.5:2) [vv] {};
      \node (r4) at (90:2) [vv] {};
      \node (r5) at (112.5:2) [vv] {};
      \node (r6) at (135:2) [vv] {};
      \node (r7) at (157.5:2) [vv] {};
      \node (r8) at (180:2) [vv] {};
      \node (r9) at (-157.5:2) [vv] {};
      \node (r10) at (-135:2) [vv] {};
      \node (r11) at (-112.5:2) [vv] {};
      \node (r12) at (-90:2) [vv] {};
      \node (r13) at (-67.5:2) [vv] {};
      \node (r14) at (-45:2) [vv] {};
      \node (r15) at (-22.5:2) [VV] {};
      \draw [R] (1) to (r); \draw [R] (r) to (r2); \draw [R] (r2) to (r3);
      \draw [R] (r3) to (r4); \draw [R] (r4) to (r5); \draw [R] (r5) to (r6);
      \draw [R] (r6) to (r7); \draw [R] (r7) to (r8); \draw [R] (r8) to (r9);
      \draw [R] (r9) to (r10); \draw [R] (r10) to (r11);
      \draw [R] (r11) to (r12); \draw [R] (r12) to (r13);
      \draw [R] (r13) to (r14); \draw [R] (r14) to (r15);
      \draw [R] (r15) to (1);
      \draw [r] (s) to (r15s);  
      \draw [bb] (1) to (s); \draw [BB] (r) to (rs);
      \draw [BB] (r2) to (r2s); \draw [BB] (r3) to (r3s);
      \draw [BB] (r4) to (r4s); \draw [BB] (r5) to (r5s);
      \draw [BB] (r6) to (r6s); \draw [BB] (r7) to (r7s);
      \draw [BB] (r8) to (r8s); \draw [BB] (r9) to (r9s);
      \draw [BB] (r10) to (r10s); \draw [BB] (r11) to (r11s);
      \draw [BB] (r12) to (r12s); \draw [BB] (r13) to (r13s);
      \draw [BB] (r14) to (r14s); \draw [bb] (r15) to (r15s);
      \node at (0:0) {$\Dih_{16}$};
    \end{scope}
    \end{tikzpicture}
    \caption{There are four semidirect products of $\Cyc_{16}$ with $\Cyc_2$, each characterized by the relation $srs=r^k$, where $k=1,7,9,15$. They result in an abelian, semidihedral, semiabelian, and dihedral group. In each case, the outer ring can be thought of as the $16^{\rm th}$ roots of unity.}\label{fig:semidirect-products}
    \end{figure}

We will start with the semidihedral group, which is the semidirect product of $\Cyc_{2^n}$ and $\Cyc_2$ defined by the automorphism $r\mapsto r^{2^{n-1}-1}$. The ``partial Cayley graph'' from Figure~\ref{fig:partial-Cayley} (middle) is completed by wiring each red arrow ``one fewer node than half way across.'' The full Cayley graphs for $\SD_8$ and $\SD_{16}$ are shown in Figure~\ref{fig:SDn-Cayley}.
  
  \tikzstyle{v-small} = [circle, draw, fill=lightgrey,inner sep=0pt, 
    minimum size=4.5mm]
  \tikzstyle{R} = [draw, very thick, Red,-stealth,bend right=10]
  \tikzstyle{R2} = [draw, very thick, Red,-stealth,bend right=15]
  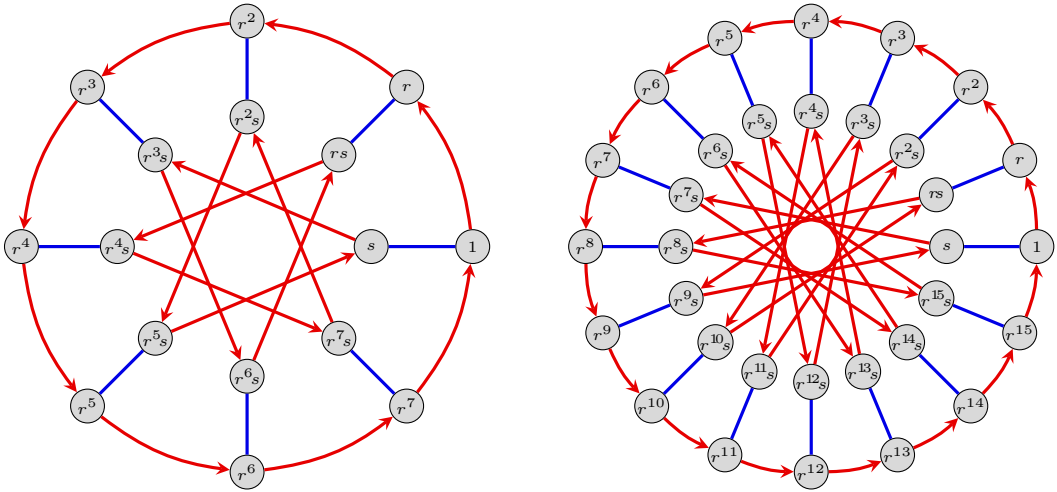
\begin{figure}[!ht]
  \begin{tikzpicture}[scale=1.5,auto]
    \tikzstyle{every node}=[font=\tiny]
    \begin{scope}[shift={(0,0)}]
      \node (s) at (0:1.1) [v-small] {$s$};
      \node (rs) at (45:1.15) [v-small] {$rs$};
      \node (r2s) at (90:1.15) [v-small] {$r^2\!s$};
      \node (r3s) at (135:1.15) [v-small] {$r^3\!s$};
      \node (r4s) at (180:1.15) [v-small] {$r^4\!s$};
      \node (r5s) at (225:1.15) [v-small] {$r^5\!s$};
      \node (r6s) at (270:1.15) [v-small] {$r^6\!s$};
      \node (r7s) at (315:1.15) [v-small] {$r^7\!s$};
      \node (1) at (0:2) [v-small] {$1$};
      \node (r) at (45:2) [v-small] {$r$};
      \node (r2) at (90:2) [v-small] {$r^2$};
      \node (r3) at (135:2) [v-small] {$r^3$};
      \node (r4) at (180:2) [v-small] {$r^4$};
      \node (r5) at (225:2) [v-small] {$r^5$};
      \node (r6) at (270:2) [v-small] {$r^6$};
      \node (r7) at (315:2) [v-small] {$r^7$};
      \draw [R2] (1) to (r);
      \draw [R2] (r) to (r2);
      \draw [R2] (r2) to (r3);
      \draw [R2] (r3) to (r4);
      \draw [R2] (r4) to (r5);
      \draw [R2] (r5) to (r6);
      \draw [R2] (r6) to (r7);
      \draw [R2] (r7) to (1);
      \draw [r] (s) to (r3s);
      \draw [r] (r3s) to (r6s);
      \draw [r] (r6s) to (rs);
      \draw [r] (rs) to (r4s);
      \draw [r] (r4s) to (r7s);
      \draw [r] (r7s) to (r2s);
      \draw [r] (r2s) to (r5s);
      \draw [r] (r5s) to (s);
      \draw [bb] (1) to (s); \draw [bb] (r) to (rs);
      \draw [bb] (r2) to (r2s); \draw [bb] (r3) to (r3s);
      \draw [bb] (r4) to (r4s); \draw [bb] (r5) to (r5s);
      \draw [bb] (r6) to (r6s); \draw [bb] (r7) to (r7s);
    \end{scope}
    \begin{scope}[shift={(5,0)}]
      \node (s) at (0:1.2) [v-small] {$s$};
      \node (rs) at (22.5:1.2) [v-small] {$r\!s$};
      \node (r2s) at (45:1.2) [v-small] {$r^2\!s$};
      \node (r3s) at (67.5:1.2) [v-small] {$r^3\!s$};
      \node (r4s) at (90:1.2) [v-small] {$r^4\!s$};
      \node (r5s) at (112.5:1.2) [v-small] {$r^5\!s$};
      \node (r6s) at (135:1.2) [v-small] {$r^6\!s$};
      \node (r7s) at (157.5:1.2) [v-small] {$r^7\!s$};
      \node (r8s) at (180:1.2) [v-small] {$r^8\!s$};
      \node (r9s) at (-157.5:1.2) [v-small] {$r^9\!s$};
      \node (r10s) at (-135:1.2) [v-small] {$r^{1\!0}\!s$};
      \node (r11s) at (-112.5:1.2) [v-small] {$r^{1\!1}\!s$};
      \node (r12s) at (-90:1.2) [v-small] {$r^{1\!2}\!s$};
      \node (r13s) at (-67.5:1.2) [v-small] {$r^{1\!3}\!s$};
      \node (r14s) at (-45:1.2) [v-small] {$r^{1\!4}\!s$};
      \node (r15s) at (-22.5:1.2) [v-small] {$r^{1\!5}\!s$};
      \node (1) at (0:2) [v-small] {$1$};
      \node (r) at (22.5:2) [v-small] {$r$};
      \node (r2) at (45:2) [v-small] {$r^2$};
      \node (r3) at (67.5:2) [v-small] {$r^3$};
      \node (r4) at (90:2) [v-small] {$r^4$};
      \node (r5) at (112.5:2) [v-small] {$r^5$};
      \node (r6) at (135:2) [v-small] {$r^6$};
      \node (r7) at (157.5:2) [v-small] {$r^7$};
      \node (r8) at (180:2) [v-small] {$r^8$};
      \node (r9) at (-157.5:2) [v-small] {$r^9$};
      \node (r10) at (-135:2) [v-small] {$r^{10}$};
      \node (r11) at (-112.5:2) [v-small] {$r^{11}$};
      \node (r12) at (-90:2) [v-small] {$r^{12}$};
      \node (r13) at (-67.5:2) [v-small] {$r^{13}$};
      \node (r14) at (-45:2) [v-small] {$r^{14}$};
      \node (r15) at (-22.5:2) [v-small] {$r^{15}$};
      \draw [R] (1) to (r); \draw [R] (r) to (r2); \draw [R] (r2) to (r3);
      \draw [R] (r3) to (r4); \draw [R] (r4) to (r5); \draw [R] (r5) to (r6);
      \draw [R] (r6) to (r7); \draw [R] (r7) to (r8); \draw [R] (r8) to (r9);
      \draw [R] (r9) to (r10); \draw [R] (r10) to (r11);
      \draw [R] (r11) to (r12); \draw [R] (r12) to (r13);
      \draw [R] (r13) to (r14); \draw [R] (r14) to (r15);
      \draw [R] (r15) to (1);
      \draw [r] (s) to (r7s); \draw [r] (r7s) to (r14s);
      \draw [r] (r14s) to (r5s); \draw [r] (r5s) to (r12s);
      \draw [r] (r12s) to (r3s); \draw [r] (r3s) to (r10s);
      \draw [r] (r10s) to (rs); \draw [r] (rs) to (r8s);
      \draw [r] (r8s) to (r15s); \draw [r] (r15s) to (r6s);
      \draw [r] (r6s) to (r13s); \draw [r] (r13s) to (r4s);
      \draw [r] (r4s) to (r11s); \draw [r] (r11s) to (r2s);
      \draw [r] (r2s) to (r9s); \draw [r] (r9s) to (s); 
      \draw [bb] (1) to (s); \draw [bb] (r) to (rs);
      \draw [bb] (r2) to (r2s); \draw [bb] (r3) to (r3s);
      \draw [bb] (r4) to (r4s); \draw [bb] (r5) to (r5s);
      \draw [bb] (r6) to (r6s); \draw [bb] (r7) to (r7s);
      \draw [bb] (r8) to (r8s); \draw [bb] (r9) to (r9s);
      \draw [bb] (r10) to (r10s); \draw [bb] (r11) to (r11s);
      \draw [bb] (r12) to (r12s); \draw [bb] (r13) to (r13s);
      \draw [bb] (r14) to (r14s); \draw [bb] (r15) to (r15s);
    \end{scope}
  \end{tikzpicture}
  \caption{The relation $srs=r^{2^{n-1}-1}$ defines the semidihedral
  group. Shown are the Cayley graphs for $\SD_8$ and $\SD_{16}$.}\label{fig:SDn-Cayley}
  \end{figure}
  
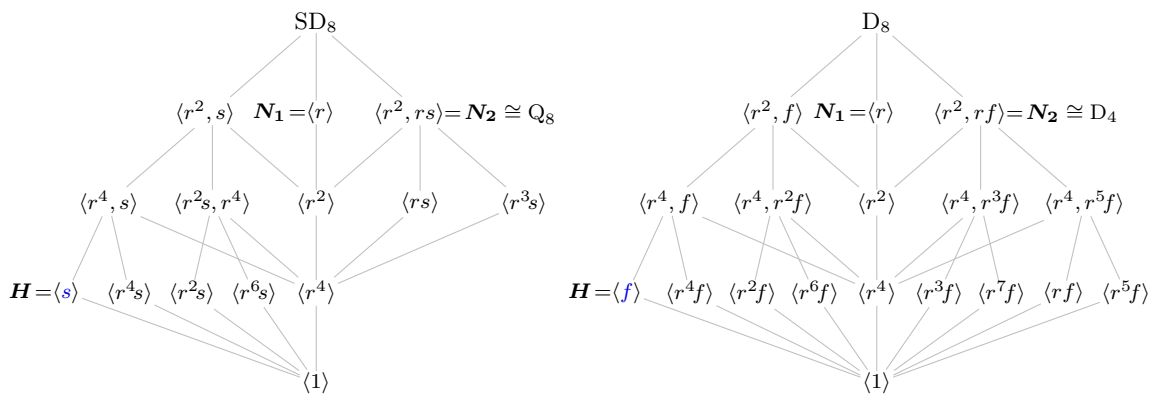
\begin{figure}[!ht]
  \begin{tikzpicture}[scale=1.2,xscale=.92,shorten >= -3pt, shorten <= -3pt]
    \begin{scope}[shift={(0,0)},scale=1]
      \tikzstyle{every node}=[font=\footnotesize]
      \node (Q16) at (0,4) {\small $\SD_8$};
      \node (r2-s) at (-1.25,3) {\!\!\!\!$\<r^2,s\>$};
      \node (r) at (0,3) {\hspace{-6mm}${\footnotesize\bm{N_1}\!=}\<r\>$};
      \node (r2-rs) at (1.25,3) {\hspace{12mm}$\<r^2,rs\>{\footnotesize=\!\bm{N_2}\cong\Quat_8}$};
      \node (s) at (-2.5,2) {$\<r^4,s\>$};
      \node (r2s) at (-1.25,2) {$\<r^2\!s,r^4\>$};
      \node (r2) at (0,2) {$\<r^2\>$};
      \node (rs) at (1.25,2) {$\<rs\>$};
      \node (r3s) at (2.5,2) {$\<r^3\!s\>$};
      \node (r4) at (0,1) {$\<r^4\>$};
      \node (a) at (-3,1) {\hspace{-6mm}${\footnotesize\bm{H}\!=}\<\Balert{s}\>$};
      \node (b) at (-2.25,1) {$\<r^4\!s\>$};
      \node (c) at (-1.5,1) {$\<r^2\!s\>$};
      \node (d) at (-.75,1) {$\<r^6\!s\>$};
      \node (1) at (0,0) {$\<1\>$};
      \draw[f] (Q16) to (r2-s); \draw[f] (Q16) to (r); \draw[f] (Q16) to (r2-rs);
      \draw[f] (r2-s) to (s); \draw[f] (r2-s) to (r2s); \draw[f] (r2-s) to (r2);
      \draw[f] (r) to (r2);
      \draw[f] (r2-rs) to (r2); \draw[f] (r2-rs) to (rs); \draw[f] (r2-rs) to (r3s);
      \draw[f] (r4) to (s); \draw[f] (r4) to (r2s); \draw[f] (r4) to (r2);
      \draw[f] (r4) to (rs); \draw[f] (r4) to (r3s); \draw[f] (1) to (r4);
      \draw[f] (a) to (1); \draw[f] (b) to (1); \draw[f] (c) to (1); \draw[f] (d) to (1);
      \draw[f] (s) to (a); \draw[f] (s) to (b); \draw[f] (r2s) to (c);
      \draw[f] (r2s) to (d); 
    \end{scope}
    \begin{scope}[shift={(6.75,0)}]
     \tikzstyle{every node}=[font=\footnotesize]
      \node (D8) at (0,4) {\small $\Dih_8$};
      \node (r2-f) at (-1.25,3) {$\<r^2,f\>$};
      \node (r) at (0,3) {\hspace{-6mm}${\footnotesize\bm{N_1}\!=}\<r\>$};
      \node (r2-rf) at (1.25,3) {\hspace{12mm}$\<r^2,rf\>{\footnotesize=\!\bm{N_2}\cong\Dih_4}$};
      \node (r4-f) at (-2.5,2) {$\<r^4,f\>$};
      \node (r4-r2f) at (-1.25,2) {$\<r^4,r^2\!f\>$};
      \node (r2) at (0,2) {$\<r^2\>$};
      \node (r4-r3f) at (1.25,2) {$\<r^4,r^3\!f\>$};
      \node (r4-r5f) at (2.5,2) {$\<r^4,r^5\!f\>$};
      \node (f) at (-3,1) {\hspace{-6mm}${\footnotesize\bm{H}\!=}\<\Balert{f}\>$};
      \node (r4f) at (-2.25,1) {$\<r^4\!f\>$};
      \node (r2f) at (-1.5,1) {$\<r^2\!f\>$};
      \node (r6f) at (-.75,1) {$\<r^6\!f\>$};
      \node (r4) at (0,1) {$\<r^4\>$};
      \node (r7f) at (1.5,1) {$\<r^7\!f\>$};
      \node (r3f) at (.75,1) {$\<r^3\!f\>$};
      \node (r5f) at (3,1) {$\<r^5\!f\>$};
      \node (rf) at (2.25,1) {$\<rf\>$};    
      \node (1) at (0,0) {$\<1\>$};
      \draw[f] (D8) to (r2-f); \draw[f] (D8) to (r); \draw[f] (D8) to (r2-rf);
      \draw[f] (r2-f) to (r4-f); \draw[f] (r2-f) to (r4-r2f); \draw[f] (r2-f) to (r2);
      \draw[f] (r2-rf) to (r4-r5f); \draw[f] (r2-rf) to (r4-r3f); \draw[f] (r2-rf) to (r2);
      \draw[f] (r2) to (r); \draw[f] (r2) to (r4);
      \draw[f] (r4-f) to (f); \draw[f] (r4-f) to (r4f); \draw[f] (r4-f) to (r4);
      \draw[f] (r4-r2f) to (r2f); \draw[f] (r4-r2f) to (r6f); \draw[f] (r4-r2f) to (r4);
      \draw[f] (r4-r3f) to (r7f); \draw[f] (r4-r3f) to (r3f); \draw[f] (r4-r3f) to (r4);
      \draw[f] (r4-r5f) to (r5f); \draw[f] (r4-r5f) to (rf); \draw[f] (r4-r5f) to (r4);
      \draw[f] (f) to (1); \draw[f] (r4f) to (1); \draw[f] (r2f) to (1); \draw[f] (r6f) to (1);
      \draw[f] (r7f) to (1); \draw[f] (r3f) to (1); \draw[f] (r5f) to (1); \draw[f] (rf) to (1);
      \draw[f] (r4) to (1);
    \end{scope}    
  \end{tikzpicture}
  \caption{The semidihedral and dihedral groups of order $16$ both have two index-$2$ subgroups that intersect trivially with an order-$2$ subgroup. Therefore, they both decompose as a semidirect product of an order-$8$ subgroup with $\Cyc_2$ in two distinct ways.} \label{fig:SD8-D8-lattice}
  \end{figure}

  Some books call the group $\SD_8$ the ``quasidihedral group,'' and denote it as $\QD_8$, or even $\QD_{16}$. We are prefer the ``semi-'' prefix for a number of reasons. First, it is a semidirect product. Second, the diquaternion group is denoted $\DQ$, and it is best to avoid mixing that with $\QD$. Third, the Cayley graph of the semidihedral group can be constructed from the dihedral group by ``rewiring'' the inner red arrow from $s\mapsto r^{-1}$ to $s\mapsto r^{2^{n-1}-1}$, which is halfway across the circle, and ``semi-'' means ``half.'' Since $r^{-1}$ and $r^{2^{n-1}-1}$ differ by a factor of $r^{2^{n-1}}=-1$, the dihedral and semidihedral groups share a number of structural properties. Their subgroup lattices, which appear in Figure~\ref{fig:SD8-D8-lattice}, are similar, with the semidihedral group missing ``the lower half'' of the non-normal order-$2$ subgroups---yet another a reason to use the ``semi-'' prefix over ``quasi-''!

The subgroup lattices of $\SD_8$ and $\Dih_8$ in Figure~\ref{fig:SD8-D8-lattice} not only look similar in structure, but they also share a number of common features. For example, the quotient by the order $2$ subgroup $\<r^4\>=\<-1\>$ is isomorphic to the dihedral group $\Dih_4$. Also, not only are they both semidirect products of $\Cyc_8$ with $\Cyc_2$, but since the ``upper-right'' index-$2$ subgroup $\<r^2,rs\>\cong \Quat_8$ (analogously, $\<r^2,rf\>\cong \Dih_4$ in $\Dih_8$) intersects one of the $\Cyc_2$ subgroups trivially, these groups are also semidirect products of this subgroup with $\Cyc_2$. 

The last semidirect product of $\Cyc_{2^n}$ with $\Cyc_2$ is defined by the relation $srs=r^{2^{n-1}+1}$. This group is so uncommon that it does not have a standard name. In a few places, it is called the \emph{modular group}, denoted by some variation of $M_2(4)$ or $M_{16}$, simply because its subgroup lattice has a certain structural property called ``modular.'' However, this name is not great, because every abelian group has a modular lattice. A few sources called it the \emph{maximum modular group}, but it is not apparent what ``maximum'' even refers to. Additionally, number theorists use the name ``modular group'' for the projective special linear group $\PSL_2(\Z)$, because of its connections to modular forms. The wonderful LMFDB group database denotes our last example $\OD_n$ for ``other dihedral.'' The name \emph{Isanowa group} also exists in the literature.

    \tikzstyle{v-small} = [circle, draw, fill=lightgrey,inner sep=0pt, 
    minimum size=5mm]
  \tikzstyle{R} = [draw, very thick, Red,-stealth,bend right=10]
  \tikzstyle{R2} = [draw, very thick, Red,-stealth,bend right=15]
  \begin{figure}[!ht]
  \begin{tikzpicture}[scale=1.5,auto]
    \tikzstyle{every node}=[font=\tiny]
    \begin{scope}[shift={(0,0)}]
      \node (s) at (0:1.1) [v-small] {$s$};
      \node (rs) at (45:1.15) [v-small] {$rs$};
      \node (r2s) at (90:1.15) [v-small] {$r^2\!s$};
      \node (r3s) at (135:1.15) [v-small] {$r^3\!s$};
      \node (r4s) at (180:1.15) [v-small] {$r^4\!s$};
      \node (r5s) at (225:1.15) [v-small] {$r^5\!s$};
      \node (r6s) at (270:1.15) [v-small] {$r^6\!s$};
      \node (r7s) at (315:1.15) [v-small] {$r^7\!s$};
      \node (1) at (0:2) [v-small] {$1$};
      \node (r) at (45:2) [v-small] {$r$};
      \node (r2) at (90:2) [v-small] {$r^2$};
      \node (r3) at (135:2) [v-small] {$r^3$};
      \node (r4) at (180:2) [v-small] {$r^4$};
      \node (r5) at (225:2) [v-small] {$r^5$};
      \node (r6) at (270:2) [v-small] {$r^6$};
      \node (r7) at (315:2) [v-small] {$r^7$};
      \draw [R2] (1) to (r);
      \draw [R2] (r) to (r2);
      \draw [R2] (r2) to (r3);
      \draw [R2] (r3) to (r4);
      \draw [R2] (r4) to (r5);
      \draw [R2] (r5) to (r6);
      \draw [R2] (r6) to (r7);
      \draw [R2] (r7) to (1);
      \draw [r] (s) to (r5s);
      \draw [r] (r5s) to (r2s);
      \draw [r] (r2s) to (r7s);
      \draw [r] (r7s) to (r4s);
      \draw [r] (r4s) to (rs);
      \draw [r] (rs) to (r6s);
      \draw [r] (r6s) to (r3s);
      \draw [r] (r3s) to (s);
      \draw [bb] (1) to (s); \draw [bb] (r) to (rs);
      \draw [bb] (r2) to (r2s); \draw [bb] (r3) to (r3s);
      \draw [bb] (r4) to (r4s); \draw [bb] (r5) to (r5s);
      \draw [bb] (r6) to (r6s); \draw [bb] (r7) to (r7s);
    \end{scope}
    \begin{scope}[shift={(5,0)}]
      \node (s) at (0:1.2) [v-small] {$s$};
      \node (rs) at (22.5:1.2) [v-small] {$r\!s$};
      \node (r2s) at (45:1.2) [v-small] {$r^2\!s$};
      \node (r3s) at (67.5:1.2) [v-small] {$r^3\!s$};
      \node (r4s) at (90:1.2) [v-small] {$r^4\!s$};
      \node (r5s) at (112.5:1.2) [v-small] {$r^5\!s$};
      \node (r6s) at (135:1.2) [v-small] {$r^6\!s$};
      \node (r7s) at (157.5:1.2) [v-small] {$r^7\!s$};
      \node (r8s) at (180:1.2) [v-small] {$r^8\!s$};
      \node (r9s) at (-157.5:1.2) [v-small] {$r^9\!s$};
      \node (r10s) at (-135:1.2) [v-small] {$r^{1\!0}\!s$};
      \node (r11s) at (-112.5:1.2) [v-small] {$r^{1\!1}\!s$};
      \node (r12s) at (-90:1.2) [v-small] {$r^{1\!2}\!s$};
      \node (r13s) at (-67.5:1.2) [v-small] {$r^{1\!3}\!s$};
      \node (r14s) at (-45:1.2) [v-small] {$r^{1\!4}\!s$};
      \node (r15s) at (-22.5:1.2) [v-small] {$r^{1\!5}\!s$};
      \node (1) at (0:2) [v-small] {$1$};
      \node (r) at (22.5:2) [v-small] {$r$};
      \node (r2) at (45:2) [v-small] {$r^2$};
      \node (r3) at (67.5:2) [v-small] {$r^3$};
      \node (r4) at (90:2) [v-small] {$r^4$};
      \node (r5) at (112.5:2) [v-small] {$r^5$};
      \node (r6) at (135:2) [v-small] {$r^6$};
      \node (r7) at (157.5:2) [v-small] {$r^7$};
      \node (r8) at (180:2) [v-small] {$r^8$};
      \node (r9) at (-157.5:2) [v-small] {$r^9$};
      \node (r10) at (-135:2) [v-small] {$r^{10}$};
      \node (r11) at (-112.5:2) [v-small] {$r^{11}$};
      \node (r12) at (-90:2) [v-small] {$r^{12}$};
      \node (r13) at (-67.5:2) [v-small] {$r^{13}$};
      \node (r14) at (-45:2) [v-small] {$r^{14}$};
      \node (r15) at (-22.5:2) [v-small] {$r^{15}$};
      \draw [R] (1) to (r); \draw [R] (r) to (r2); \draw [R] (r2) to (r3);
      \draw [R] (r3) to (r4); \draw [R] (r4) to (r5); \draw [R] (r5) to (r6);
      \draw [R] (r6) to (r7); \draw [R] (r7) to (r8); \draw [R] (r8) to (r9);
      \draw [R] (r9) to (r10); \draw [R] (r10) to (r11);
      \draw [R] (r11) to (r12); \draw [R] (r12) to (r13);
      \draw [R] (r13) to (r14); \draw [R] (r14) to (r15);
      \draw [R] (r15) to (1);
      \draw [r] (s) to (r9s); \draw [r] (r9s) to (r2s);
      \draw [r] (r2s) to (r11s); \draw [r] (r11s) to (r4s);
      \draw [r] (r4s) to (r13s); \draw [r] (r13s) to (r6s);
      \draw [r] (r6s) to (r15s); \draw [r] (r15s) to (r8s);
      \draw [r] (r8s) to (rs); \draw [r] (rs) to (r10s);
      \draw [r] (r10s) to (r3s); \draw [r] (r3s) to (r12s);
      \draw [r] (r12s) to (r5s); \draw [r] (r5s) to (r14s);
      \draw [r] (r14s) to (r7s); \draw [r] (r7s) to (s); 
      \draw [bb] (1) to (s); \draw [bb] (r) to (rs);
      \draw [bb] (r2) to (r2s); \draw [bb] (r3) to (r3s);
      \draw [bb] (r4) to (r4s); \draw [bb] (r5) to (r5s);
      \draw [bb] (r6) to (r6s); \draw [bb] (r7) to (r7s);
      \draw [bb] (r8) to (r8s); \draw [bb] (r9) to (r9s);
      \draw [bb] (r10) to (r10s); \draw [bb] (r11) to (r11s);
      \draw [bb] (r12) to (r12s); \draw [bb] (r13) to (r13s);
      \draw [bb] (r14) to (r14s); \draw [bb] (r15) to (r15s);
    \end{scope}
  \end{tikzpicture}
  \caption{Cayley graphs of the semiabelian groups $\SA_8$ and $\SA_{16}$.}\label{fig:SAn-Cayley}
  \end{figure}
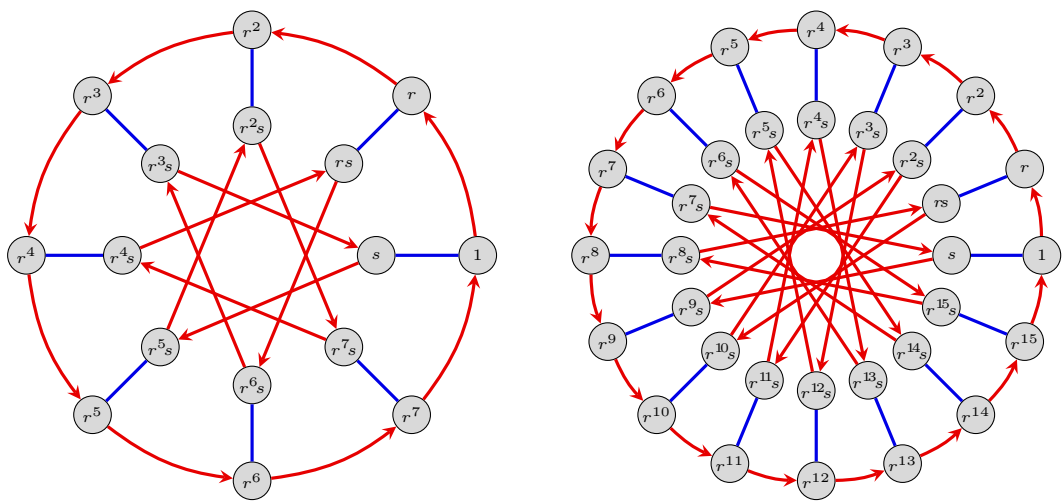

In Fall 2021, a student from the author's undergraduate abstract algebra class suggested the name ``semiabelian'' for this group. This is a fantastic term for a number of reasons. For one, just like how the presentations of the dihedral and semidihedral groups differ by $r^{2^{n-1}}=-1$, a semicircle around the inner ring of nodes (recall Figure~\ref{fig:semidirect-products}), so do the abelian and semiabelian groups. Cayley graphs of $\SA_8$ and $\SA_{16}$ are shown in Figure~\ref{fig:SAn-Cayley}. A quick Google search shows that the term ``semiabeilan'' has been used independently by (at least) two authors in isolated research papers, but in both instances, that term is general enough that it encompasses the groups $\SA_n$ as special cases. Our name works!

\begin{figure}[!ht]
  \tikzstyle{v-tiny} = [circle, draw, fill=lightgrey,inner sep=0pt, 
    minimum size=4.75mm]
  \begin{tikzpicture}[scale=.85]
   \tikzstyle{every node}=[font=\tiny]
    \begin{scope}[shift={(0,0)}]
      \node at (3.5,3.1) {\small $\Cyc_8\!\times\!\Cyc_2$}; 
      \node (10) at (0,2.5) [v-tiny] {$1$};
      \node (11) at (1,2.5) [v-tiny] {$r$};
      \node (12) at (2,2.5) [v-tiny] {$r^2$};
      \node (13) at (3,2.5) [v-tiny] {$r^3$};
      \node (14) at (4,2.5) [v-tiny] {$r^4$};
      \node (15) at (5,2.5) [v-tiny] {$r^5$};
      \node (16) at (6,2.5) [v-tiny] {$r^6$};
      \node (17) at (7,2.5) [v-tiny] {$r^7$};
      \node (00) at (0,0) [v-tiny] {$s$};
      \node (01) at (1,0) [v-tiny] {$sr$};
      \node (02) at (2,0) [v-tiny] {$sr^2$};
      \node (03) at (3,0) [v-tiny] {$sr^3$};
      \node (04) at (4,0) [v-tiny] {$sr^4$};
      \node (05) at (5,0) [v-tiny] {$sr^5$};
      \node (06) at (6,0) [v-tiny] {$sr^6$};
      \node (07) at (7,0) [v-tiny] {$sr^7$};
      \draw [r] (00) to (01); \draw [r] (01) to (02); \draw [r] (02) to (03);
      \draw [r] (03) to (04); \draw [r] (04) to (05); \draw [r] (05) to (06);
      \draw [r] (06) to (07); \draw [r] (07) to [bend left=27] (00);
      \draw [r] (10) to (11); \draw [r] (11) to (12); \draw [r] (12) to (13);
      \draw [r] (13) to (14); \draw [r] (14) to (15); \draw [r] (15) to (16);
      \draw [r] (16) to (17); \draw [r] (17) to [bend right=27] (10);
      \draw [bb] (10) to (00);
      \draw [bb] (11) to (01);
      \draw [bb] (12) to (02);
      \draw [bb] (13) to (03);
      \draw [bb] (14) to (04);
      \draw [bb] (15) to (05);
      \draw [bb] (16) to (06);
      \draw [bb] (17) to (07);
    \end{scope}
    \begin{scope}[shift={(9,0)}]
      \node at (3.5,3.1) {\small $\Dih_8$}; 
      \node (10) at (0,2.5) [v-tiny] {$1$};
      \node (11) at (1,2.5) [v-tiny] {$r$};
      \node (12) at (2,2.5) [v-tiny] {$r^2$};
      \node (13) at (3,2.5) [v-tiny] {$r^3$};
      \node (14) at (4,2.5) [v-tiny] {$r^4$};
      \node (15) at (5,2.5) [v-tiny] {$r^5$};
      \node (16) at (6,2.5) [v-tiny] {$r^6$};
      \node (17) at (7,2.5) [v-tiny] {$r^7$};
      \node (00) at (0,0) [v-tiny] {$s$};
      \node (01) at (1,0) [v-tiny] {$sr$};
      \node (02) at (2,0) [v-tiny] {$sr^2$};
      \node (03) at (3,0) [v-tiny] {$sr^3$};
      \node (04) at (4,0) [v-tiny] {$sr^4$};
      \node (05) at (5,0) [v-tiny] {$sr^5$};
      \node (06) at (6,0) [v-tiny] {$sr^6$};
      \node (07) at (7,0) [v-tiny] {$sr^7$};
      \draw [r] (00) to (01); \draw [r] (01) to (02); \draw [r] (02) to (03);
      \draw [r] (03) to (04); \draw [r] (04) to (05); \draw [r] (05) to (06);
      \draw [r] (06) to (07); \draw [r] (07) to [bend left=27] (00);
      \draw [r] (10) to (11); \draw [r] (11) to (12); \draw [r] (12) to (13);
      \draw [r] (13) to (14); \draw [r] (14) to (15); \draw [r] (15) to (16);
      \draw [r] (16) to (17); \draw [r] (17) to [bend right=27] (10);
      \draw [bb] (10) to (00);
      \draw [bb] (11) to (07);
      \draw [bb] (12) to (06);
      \draw [bb] (13) to (05);
      \draw [bb] (14) to (04);
      \draw [bb] (15) to (03);
      \draw [bb] (16) to (02);
      \draw [bb] (17) to (01);
    \end{scope}
    \begin{scope}[shift={(0,-5)}]
      \node at (3.5,3.1) {\small $\SA_8$}; 
      \node (10) at (0,2.5) [v-tiny] {$1$};
      \node (11) at (1,2.5) [v-tiny] {$r$};
      \node (12) at (2,2.5) [v-tiny] {$r^2$};
      \node (13) at (3,2.5) [v-tiny] {$r^3$};
      \node (14) at (4,2.5) [v-tiny] {$r^4$};
      \node (15) at (5,2.5) [v-tiny] {$r^5$};
      \node (16) at (6,2.5) [v-tiny] {$r^6$};
      \node (17) at (7,2.5) [v-tiny] {$r^7$};
      \node (00) at (0,0) [v-tiny] {$s$};
      \node (01) at (1,0) [v-tiny] {$sr$};
      \node (02) at (2,0) [v-tiny] {$sr^2$};
      \node (03) at (3,0) [v-tiny] {$sr^3$};
      \node (04) at (4,0) [v-tiny] {$sr^4$};
      \node (05) at (5,0) [v-tiny] {$sr^5$};
      \node (06) at (6,0) [v-tiny] {$sr^6$};
      \node (07) at (7,0) [v-tiny] {$sr^7$};
      \draw [r] (00) to (01); \draw [r] (01) to (02); \draw [r] (02) to (03);
      \draw [r] (03) to (04); \draw [r] (04) to (05); \draw [r] (05) to (06);
      \draw [r] (06) to (07); \draw [r] (07) to [bend left=27] (00);
      \draw [r] (10) to (11); \draw [r] (11) to (12); \draw [r] (12) to (13);
      \draw [r] (13) to (14); \draw [r] (14) to (15); \draw [r] (15) to (16);
      \draw [r] (16) to (17); \draw [r] (17) to [bend right=27] (10);
      \draw [bb] (10) to (00);
      \draw [bb] (11) to (05);
      \draw [bb] (12) to (02);
      \draw [bb] (13) to (07);
      \draw [bb] (14) to (04);
      \draw [bb] (15) to (01);
      \draw [bb] (16) to (06);
      \draw [bb] (17) to (03);
    \end{scope}
    \begin{scope}[shift={(9,-5)}]
      \node at (3.5,3.1) {\small $\SD_8$}; 
      \node (10) at (0,2.5) [v-tiny] {$1$};
      \node (11) at (1,2.5) [v-tiny] {$r$};
      \node (12) at (2,2.5) [v-tiny] {$r^2$};
      \node (13) at (3,2.5) [v-tiny] {$r^3$};
      \node (14) at (4,2.5) [v-tiny] {$r^4$};
      \node (15) at (5,2.5) [v-tiny] {$r^5$};
      \node (16) at (6,2.5) [v-tiny] {$r^6$};
      \node (17) at (7,2.5) [v-tiny] {$r^7$};
      \node (00) at (0,0) [v-tiny] {$s$};
      \node (01) at (1,0) [v-tiny] {$sr$};
      \node (02) at (2,0) [v-tiny] {$sr^2$};
      \node (03) at (3,0) [v-tiny] {$sr^3$};
      \node (04) at (4,0) [v-tiny] {$sr^4$};
      \node (05) at (5,0) [v-tiny] {$sr^5$};
      \node (06) at (6,0) [v-tiny] {$sr^6$};
      \node (07) at (7,0) [v-tiny] {$sr^7$};
      \draw [r] (00) to (01); \draw [r] (01) to (02); \draw [r] (02) to (03);
      \draw [r] (03) to (04); \draw [r] (04) to (05); \draw [r] (05) to (06);
      \draw [r] (06) to (07); \draw [r] (07) to [bend left=27] (00);
      \draw [r] (10) to (11); \draw [r] (11) to (12); \draw [r] (12) to (13);
      \draw [r] (13) to (14); \draw [r] (14) to (15); \draw [r] (15) to (16);
      \draw [r] (16) to (17); \draw [r] (17) to [bend right=27] (10);
      \draw [bb] (10) to (00);
      \draw [bb] (11) to (03);
      \draw [bb] (12) to (06);
      \draw [bb] (13) to (01);
      \draw [bb] (14) to (04);
      \draw [bb] (15) to (07);
      \draw [bb] (16) to (02);
      \draw [bb] (17) to (05);
    \end{scope}
\end{tikzpicture}
  \caption{Another layout of the four semidirect products of $\Cyc_8$ with $\Cyc_2$. Notice how the abelian and semiabelian groups differ by only four bidirected $s$-edges, as do the dihedral and semidihedral groups.}\label{fig:C8xC2-semidirect}
\end{figure}
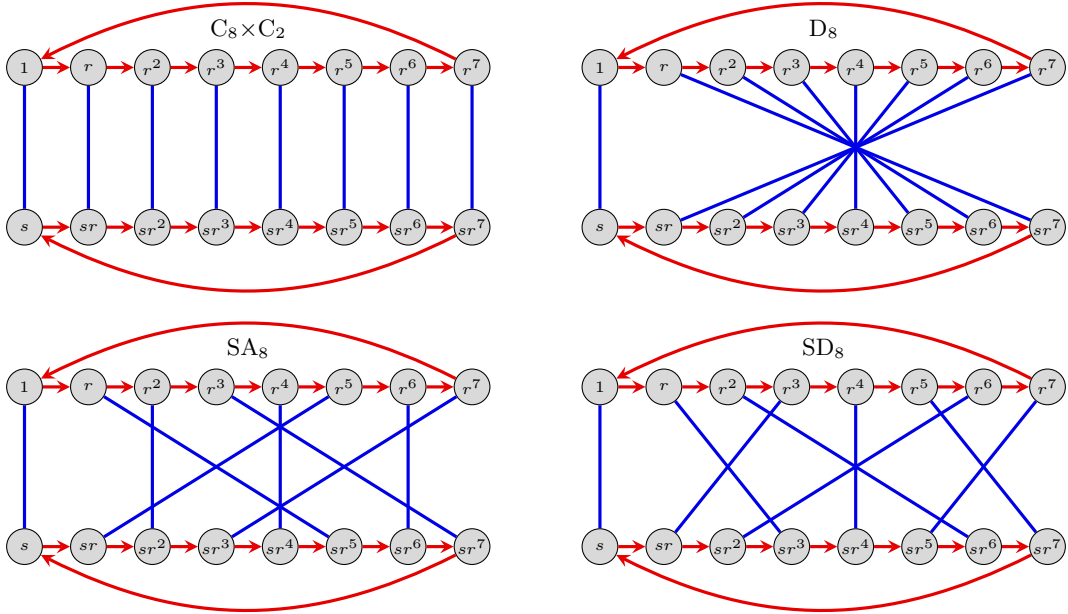

Despite the apparent structural similarities of the Cayley graphs of the semiabelian and semidihedral groups---the only difference is that the inner red edges are reversed, the semiabelian groups are more closely related to the abelian groups $\Cyc_{2^n}\times \Cyc_2$. One way to see this is to lay out their Cayley graphs differently, so the two cosets of $\<r\>$ are parallel. This is done for the groups of order $16$ in Figure~\ref{fig:C8xC2-semidirect}. Notice how the abelian and semiabelian groups differ in four of the eight bidirected $s$-edges, as do the dihedral and semidihedral groups. In contrast, the other four pairs that can be formed from these four groups differ in six of the eight $s$-edges---only the ones adjacent to $1$ and $r^4=-1$ are the same, and these are forced because these elements are central. 

   Beyond just the Cayley graphs, it is insightful to look at the representations of these groups by $2\times 2$ complex matrices, and specifically, ones that contain the reflection matrix $F$ for the order-$2$ generator. With this choice, the abelian and semiabelian groups are isomorphic to the following groups of matrices:
  \[
  \Cyc_n\times \Cyc_2\cong\left\<
  \begin{bmatrix}\zeta_n & 0 \\ 0 & \zeta_n\end{bmatrix},\;
  \begin{bmatrix}0 & 1 \\ 1 & 0\end{bmatrix}\right\>,
  \qquad
  \SA_n\cong\left\<\begin{bmatrix}\zeta_n & 0 \\ 0 & -\zeta_n\end{bmatrix},\;
  \begin{bmatrix}0 & 1 \\ 1 & 0\end{bmatrix}\right\>.
  \]
  In other words, they differ only by a negative sign in the $(2,2)$-entry of the matrix for $r$, which incorporates how the last relation in their presentations differs by $r^{2^{n-1}}=-1$. Similarly, the representations of the dihedral and semidihedral groups also differ by a negative sign in the $(2,2)$-entry:
   \[  
    \Dih_n\cong\left\<
  \begin{bmatrix}\zeta_n & 0 \\ 0 & \overline{\zeta}_n\end{bmatrix},\;
  \begin{bmatrix}0 & 1 \\ 1 & 0\end{bmatrix}\right\>,\qquad
  \SD_n\cong\left\<\begin{bmatrix}\zeta_n & 0 \\ 0 & -\overline{\zeta}_n\end{bmatrix},\; 
  \begin{bmatrix}0 & 1 \\ 1 & 0\end{bmatrix}\right\>.
  \]
  
  We will conclude this section with examples of the cycle graphs of our four semidirect products. Those for the three nonabelian groups, $\Dih_8$, $\SD_8$, and $\SA_8$, are shown in Figure~\ref{fig:D8-SD8-SA8-cycle-graphs}. The nodes are colored to track certain subsets that differ across these groups.
  Notice how all three groups have the same $8$-cycle, but they differ in the number of elements of order $2$. The group $\SA_8$ has two cyclic subgroups of order $8$. It is also worth comparing these cycle graphs to that of the generalized quaternion group $\Quat_{16}=\Dic_8$ in Figure~\ref{fig:Q16-DQ8-cycle-graphs}.

\tikzstyle{v-small} = [circle, draw, fill=lightgrey,inner sep=0pt,
    minimum size=4.75mm]
  \tikzstyle{v-r} = [circle, draw, fill=lred,inner sep=0pt, minimum size=4.75mm]
  \tikzstyle{v-g} = [circle, draw, fill=lightgrey,inner sep=0pt, minimum size=4.75mm]
  \tikzstyle{v-y} = [circle, draw, fill=lyellow,inner sep=0pt, minimum size=4.75mm]
    \tikzstyle{v-b} = [circle, draw, fill=lblue,inner sep=0pt, minimum size=4.75mm]
\tikzstyle{v-p} = [circle, draw, fill=conj-purple,inner sep=0pt, minimum size=4.75mm]
\tikzstyle{edge} = [draw,thick]
\tikzstyle{cy2} = [draw,thick]
\tikzstyle{cy2-r} = [draw,thick]
\tikzstyle{cy2-rs} = [draw,Blue,thick]
\tikzstyle{cy2-r2s} = [draw,Purple,thick]
\tikzstyle{cy2-r3s} = [draw,darkgreen,thick]
  \begin{figure}[!ht]
  \begin{tikzpicture}[scale=.6,auto]
    \tikzstyle{every node}=[font=\scriptsize]
    \begin{scope}[shift={(0,0)}]
      \node (1) at (0,3) [v-r] {$1$};
      \node (r) at (.75,1.5) [v-r] {$r$};
      \node (r2) at (1.5,0) [v-r] {$r^2$};
      \node (r3) at (.75,-1.5) [v-r] {$r^3$};
      \node (r4) at (0,-3) [v-r] {$r^4$};
      \node (r5) at (-.75,-1.5) [v-r] {$r^5$};
      \node (r6) at (-1.5,0) [v-r] {$r^6$};
      \node (r7) at (-.75,1.5) [v-r] {$r^7$};
      \node (s) at (.5,5) [v-b] {$s$};
      \node (rs) at (1.5,5) [v-y] {$rs$};
      \node (r2s) at (2.5,5) [v-p] {$r^2\!s$};
      \node (r3s) at (3.5,5) [v-y] {$r^3\!s$};
      \node (r4s) at (-.5,5) [v-b] {$r^4\!s$};
      \node (r5s) at (-1.5,5) [v-y] {$r^5\!s$};
      \node (r6s) at (-2.5,5) [v-p] {$r^6\!s$};
      \node (r7s) at (-3.5,5) [v-y] {$r^7\!s$};
      \draw [cy2] (1) to (r); \draw [cy2] (r) to (r2);
      \draw [cy2] (r2) to (r3); \draw [cy2] (r3) to (r4);
      \draw [cy2] (r4) to (r5); \draw [cy2] (r5) to (r6);
      \draw [cy2] (r6) to (r7); \draw [cy2] (r7) to (1);
      \draw [cy2] (1) to (s); \draw [cy2] (1) to (r2s);
      \draw [cy2] (1) to (r4s); \draw [cy2] (1) to (r6s);
      \draw [cy2] (1) to (rs); \draw [edge, bend right=4] (1) to (r3s);
      \draw [cy2] (1) to (r5s); \draw [edge,bend left=4] (1) to (r7s);
      \node at (0,0) {\normalsize $\Dih_8$};
    \end{scope}
    \begin{scope}[shift={(8.25,0)}]
      \node (1) at (0,3) [v-r] {$1$};
      \node (r) at (.75,1.5) [v-r] {$r$};
      \node (r2) at (1.5,0) [v-r] {$r^2$};
      \node (r3) at (.75,-1.5) [v-r] {$r^3$};
      \node (r4) at (0,-3) [v-r] {$r^4$};
      \node (r5) at (-.75,-1.5) [v-r] {$r^5$};
      \node (r6) at (-1.5,0) [v-r] {$r^6$};
      \node (r7) at (-.75,1.5) [v-r] {$r^7$};
      \node (rs) at (2.5,0) [v-y] {$rs$};
      \node (r3s) at (3.5,0) [v-y] {$r^3\!s$};
      \node (r7s) at (-2.5,0) [v-y] {$r^7\!s$};
      \node (r5s) at (-3.5,0) [v-y] {$r^5\!s$};
      \node (s) at (1,5) [v-b] {$s$};
      \node (r2s) at (3,5) [v-p] {$r^2\!s$};
      \node (r4s) at (-1,5) [v-b] {$r^4\!s$};
      \node (r6s) at (-3,5) [v-p] {$r^6\!s$};
      \draw [cy2-r] (1) to (r); \draw [cy2-r] (r) to (r2);
      \draw [cy2-r] (r2) to (r3); \draw [cy2-r] (r3) to (r4);
      \draw [cy2-r] (r4) to (r5); \draw [cy2-r] (r5) to (r6);
      \draw [cy2-r] (r6) to (r7); \draw [cy2-r] (r7) to (1);
      \draw [cy2] (1) to (s); \draw [cy2] (1) to (r2s);
      \draw [cy2] (1) to (r4s); \draw [cy2] (1) to (r6s);
      \draw [cy2-rs, bend left=8] (1) to (rs);
      \draw [cy2-rs, bend right=8] (r4) to (rs);
      \draw [cy2-r3s, bend left=5] (1) to (r3s);
      \draw [cy2-r3s, bend right=5] (r4) to (r3s);
      \draw [cy2-r3s, bend right=5] (1) to (r5s);
      \draw [cy2-r3s, bend left=5] (r4) to (r5s);
      \draw [cy2-rs, bend right=8] (1) to (r7s);
      \draw [cy2-rs, bend left=8] (r4) to (r7s);      
      \node at (0,0) {\normalsize $\SD_8$};
    \end{scope}
    \begin{scope}[shift={(16.5,0)}]
      \node (1) at (0,3) [v-r] {$1$};
      \node (r) at (1.5,1.5) [v-r] {$r$};
      \node (r2) at (3,0) [v-r] {$r^2$};
      \node (r3) at (1.5,-1.5) [v-r] {$r^3$};
      \node (r4) at (0,-3) [v-r] {$r^4$};
      \node (r5) at (-1.5,-1.5) [v-r] {$r^5$};
      \node (r6) at (-3,0) [v-r] {$r^6$};
      \node (r7) at (-1.5,1.5) [v-r] {$r^7$};
      \node (rs) at (-3,3) [v-y] {$rs$};
      \node (r3s) at (-3,-3) [v-y] {$r^3\!s$};
      \node (r7s) at (3,3) [v-y] {$r^7\!s$};
      \node (r5s) at (3,-3) [v-y] {$r^5\!s$};
      \node (r2s) at (1.25,0) [v-p] {$r^2\!s$};
      \node (r6s) at (-1.25,0) [v-p] {$r^6\!s$};
      \node (s) at (1.5,5) [v-b] {$s$};
      \node (r4s) at (-1.5,5) [v-b] {$r^4\!s$};

      \draw [cy2] (1) to (r); \draw [cy2] (r) to (r2);
      \draw [cy2] (r2) to (r3); \draw [cy2] (r3) to (r4);
      \draw [cy2] (r4) to (r5); \draw [cy2] (r5) to (r6);
      \draw [cy2] (r6) to (r7); \draw [cy2] (r7) to (1);
      \draw [cy2] (1) to (s); \draw [cy2] (1) to (r4s);
      \draw [cy2-r2s] (1) to (r6s); \draw [cy2-r2s] (r6s) to (r4);
      \draw [cy2-r2s] (1) to (r2s); \draw [cy2-r2s] (r2s) to (r4);
      \draw [cy2-rs] (1) to (r7s); \draw [cy2-rs] (r7s) to (r2);
      \draw [cy2-rs] (r2) to (r5s); \draw [cy2-rs] (r5s) to (r4);
      \draw [cy2-rs] (r4) to (r3s); \draw [cy2-rs] (r3s) to (r6);
      \draw [cy2-rs] (r6) to (rs); \draw [cy2-rs] (rs) to (1);
      \node at (0,0) {\normalsize $\SA_8$};
    \end{scope}    
\end{tikzpicture}
\caption{Cycle graphs of the three nonabelian semidirect products of $\Cyc_8$ with $\Cyc_2$. Shaded and colored nodes have the same orbit structure in all three groups.
}\label{fig:D8-SD8-SA8-cycle-graphs}
  \end{figure}
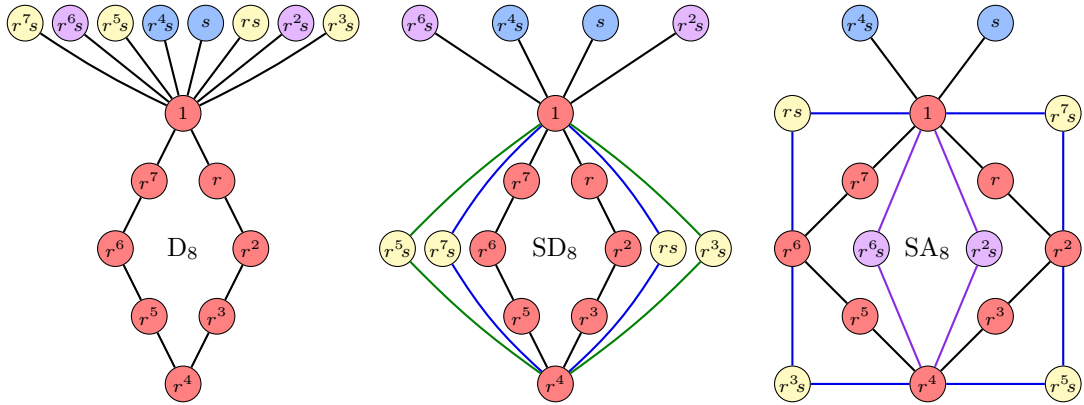
    
 At last, we need to address the abelian group $\Cyc_8\times \Cyc_2$, whose cycle graph seems to be conspicuously missing from Figure~\ref{fig:D8-SD8-SA8-cycle-graphs}. Actually, this is not the case, because it has the exact same cycle graph as $\SA_8$! This means that it has two cyclic subgroups of order $8$, $\<r\>$ and $\<rs\>$, and a third order-$8$ subgroup isomorphic to $\Cyc_4\times \Cyc_2$. This is enough to force these two groups to have identical subgroup lattices. We have seen this lattice---it is the mystery one that appeared on the first page this article, with $7$ unicorns among its $11$ subgroups, and either $9$ or $11$ normal subgroups. Alas, we have reached the solution to our puzzle. This is the smallest instance of two nonisomorphic groups that have the same subgroup lattice.\footnote{By ``the same,'' we mean that the subgroup lattices, where each edge $H\leq K$ is weighted by the index $[H:K]$, are identical. In particular, this means that $\Cyc_2$ and $\Cyc_3$ have different subgroup lattices.} It is worth noting that in general, having the same cycle graph is not sufficient for having the same subgroup lattice. We encourage the reader to verify, perhaps with the help of LMFDB, that the abelian group $\Cyc_4\times \Cyc_2^2$ has the same cycle graph as $\DQ_8$, but a different subgroup lattice.

\section{Concluding remarks}

We hope that readers take away several key messages from this article. First, there are many interesting and accessible finite groups beyond just $\Cyc_n$, $\Dih_n$, $\Sym_n$, $\Alt_n$, and $\Quat_8$, which make up the majority of examples in a standard introductory algebra class. The groups we have just seen exhibit a number of interesting properties that the aforementioned ones do not, which are illuminating to show in class. For example, the semidihedral group $\SD_8$ is the Galois group of $f(x)=x^8-2$, and provides an even richer ``Goldilocks example'' of a concrete Galois group than the classic example of $\Gal(x^4-2)\cong \Dih_4$. The semiabelian groups give an explicit example of how nonisomorphic groups can have an identical subgroup structure. The diquaternion groups provide a smallish example of an interesting central product, and an example of a group that decomposes into a semidirect product three different ways.\footnote{The dihedral group $\Dih_6$ also decomposes as a semidirect product three ways, \emph{and} as a direct product.} In contrast, the generalized quaternion groups are examples of groups that do \emph{not} break up as semidirect products at all. The dicyclic groups arise from a natural construction, that many algebra students and instructors (including the author of this article), had simply never considered. 

Another take-away is that there is no way we could have communicated these ideas as effectively without the visual tools. It is unfortunate that these are largely absent from classrooms and books, with a few exceptions, like \cite{carter2009visual,ernst2023inquiry,macauley2023visual}. This poses questions about traditional abstract algebra pedagogy, and its effectiveness in communicating key concepts. As a beginning graduate student said to the author in an introductory email, seeing these visuals is like having a superpower, and there is ``no going back.'' Furthermore, the ideas in this article are just the tip of the iceberg in terms of the mileage that can be achieved in an undergraduate or even graduate abstract algebra class by supplementing with visuals. Instructors should take full advantage of resources such as the GroupNames and LMFDB websites, and incorporate them into assignments and assessments. The LMFDB has a search capability with many parameters. In a split-second, it will return a complete list of all (up to a reasonable size) non-nilpotent solvable groups with derived length $3$, center $Z(G)\cong \Cyc_4$, that decompose as a semidirect product but not as a direct product. It also displays the (reduced) subgroup lattices for those that are not too big, and provides additional information about the groups and subgroups. Students can explore questions like what happens if one tries to ``dihedralize the dicyclic groups,'' or what matrices result from the representations of $\SD_n$ and $\SA_n$ if $n$ is not a power of $2$. The possibilities are endless.

We will conclude with a parting thought about the groups discussed in this article. For any $n\geq 4$, there are exactly six groups of order $2^n$ that have an index-$2$ cyclic subgroup: the abelian groups $\Cyc_{2^n}$ and $\Cyc_{2^{n-1}}\times \Cyc_2$, the three nonabelian semidirect products (dihedral, semidihedral, and semiabelian), and the generalized quaternion group $\Quat_{2^n}$ (see Chapter 5 of \cite{gorenstein2007finite} for a proof). A fun fact that usually comes as quite a surprise to even research algebraists is just how many groups there are of order $2^n$. For example, there are 49,487,367,289 groups of order $1024$,\footnote{In a 2002 paper \cite{besche2002millennium}, this number was reported to be 49,487,365,422, but a correction appeared in 2022 \cite{burrell2021number}.} which is over $99\%$ of all groups of order at most $2000$. Yet, only six of these have an index-$2$ subgroup (equivalently, an element of order $512$). There are still just six among the $1.77\times 10^{15}$ groups of order 2048,\footnote{John Conway said that the human race will never know the exact number of groups of order $2048$, but in \cite{conway2008counting} he writes that it exceeds the $1,\!774,\!274,\!116,\!992,\!170$ exponent-2 class 2 groups, and that the true number shares the first three digits.} and in this paper, we learned how to construct all of them. As far as we know, it is still an open problem to prove that ``almost all'' finite groups are $2$-groups, in the asymptotic sense, despite this being an absolute certainty.

\end{document}